\RequirePackage{fix-cm}

\documentclass[a4]{article} 

\usepackage{multirow}

\usepackage{graphicx}
\usepackage{amssymb}
\usepackage{amsmath}
\usepackage{amsfonts}
\usepackage{url}
\usepackage[pdftex]{thumbpdf}
\usepackage[pdftex,
pdfstartview=FitH,bookmarks=true,bookmarksnumbered=true,
bookmarksopen=true,hypertexnames=false,breaklinks=true,
colorlinks=true,linkcolor=blue,anchorcolor=blue,
citecolor=blue,filecolor=blue,menucolor=blue]{hyperref}%
\newtheorem{assumption}{Assumption}
\newtheorem{algorithm}{Algorithm}
\newtheorem{corollary}{Corollary}
\newtheorem{lemma}{Lemma}
\newtheorem{remark}{Remark}
\newtheorem{proposition}{Proposition}

\begin{document}

\title{{A}daptive-{M}ultilevel {BDDC} and its parallel implementation}

\author{
Bed\v{r}ich~Soused\'{\i}k\footnote{
      Department of Aerospace and Mechanical Engineering, University of Southern California, 
      Los Angeles, CA 90089-2531, USA,
      and Institute of Thermomechanics, Academy of Sciences of the Czech Republic,
      Dolej\v{s}kova 1402/5, CZ - 182~00 Prague~8, Czech Republic.
      sousedik@usc.edu 
}
\and
Jakub~\v{S}\'{\i}stek\footnote{
      Institute of Mathematics, Academy of Sciences of the Czech Republic, 
      \v Zitn\' a 25, 115 67 Prague 1, Czech Republic. 
      Corresponding author, Tel.: +420 222 090 710; Fax: +420 222 211 638. 
      sistek@math.cas.cz
}
\and 
Jan~Mandel\footnote{
      Department of Mathematical and Statistical Sciences, University of Colorado Denver, 
      Denver, CO 80217-3364, USA.
      jan.mandel@ucdenver.edu 
}
}

\date{}

\maketitle

\begin{abstract}
We combine the adaptive and multilevel approaches to the BDDC and formulate a method
which allows an adaptive selection of constraints on each decomposition level.
We also present a strategy for the solution of local eigenvalue problems in the adaptive algorithm using
the LOBPCG method with a preconditioner based on standard components of the BDDC.
The effectiveness of the method is illustrated on several engineering problems.
It appears that the {A}daptive-{M}ultilevel {BDDC} algorithm is able to effectively detect troublesome parts
on each decomposition level and improve convergence of the method. 
The developed open-source parallel implementation shows a good scalability as well as applicability to very large problems and core counts.

\vskip 5mm
\centerline{
Dedicated to Professor Ivo Marek on the occasion of his 80th birthday
}
\end{abstract}

\section{Introduction}

The \emph{Balancing Domain Decomposition by Constraints} (BDDC) was
developed by Dohrmann~\cite{Dohrmann-2003-PSC} as a~primal alternative to
the \emph{Finite Element Tearing and Interconnecting - Dual, Primal}
(FETI-DP) by Farhat et al.~\cite{Farhat-2000-SDP}. Both methods use
constraints to impose equality of new `coarse' variables on substructure interfaces, such as
values at substructure corners or weighted averages over edges and faces.
Primal variants of the FETI-DP were also independently proposed by Cros~\cite%
{Cros-2003-PSC} and by Fragakis and Papadrakakis~\cite{Fragakis-2003-MHP}.
It has been shown in~\cite{Mandel-2007-BFM,Sousedik-2008-EPD} that these
methods are in fact the same as BDDC. Polylogarithmic condition number
bounds for FETI-DP were first proved in~\cite{Mandel-2001-CDP} and
generalized to the case of coefficient jumps between substructures in \cite%
{Klawonn-2002-DPF}. The same bounds were obtained for BDDC in~\cite%
{Mandel-2003-CBD,Mandel-2005-ATP}. A proof that the eigenvalues of the
preconditioned operators of both methods are actually the same except for
the eigenvalues equal to one was given in~\cite{Mandel-2005-ATP} and then
simplified in~\cite{Brenner-2007-BFW,Li-2006-FBB,Mandel-2007-BFM}. FETI-DP,
and, equivalently,\ BDDC are quite robust. It can be proved that the
condition number remains bounded even for large classes of subdomains with
rough interfaces in 2D \cite{Klawonn-2008-AFA,Widlund-2008-AIS} as well as
in many cases of strong discontinuities of coefficients, including some
configurations when the discontinuities cross substructure boundaries \cite%
{Pechstein-2008-AFM,Pechstein-2009-AFM}. However, the condition number does
deteriorate in many situations of practical importance and an adaptive
method is warranted.

Adaptive enrichment for BDDC and FETI-DP was proposed in \cite%
{Mandel-2006-ACS,Mandel-2007-ASF}, with the added coarse functions built
from eigenproblems based on adjacent pairs of substructures in 2D formulated
in terms of FETI-DP\ operators. The algorithm has been developed directly in
terms of BDDC operators and extended to 3D by Mandel, Soused\'{\i}k and \v{S}%
\'{\i}stek \cite{Mandel-2012-ABT,Sousedik-2008-CDD}, resulting in a much
simplified formulation and implementation with global matrices, no explicit
coarse problem, and getting much of its parallelism through the direct
solver used to solve an auxiliary decoupled system. The only requirement for
all these versions of the adaptive algorithms is that there is a sufficient
number of corner constraints to prevent rigid body motions between any pair
of adjacent substructures. This requirement has been recognized in
other contexts~\cite{Broz-2009-ACN,Lesoinne-2003-FCS}, and in the context of
BDDC\ by Dohrmann~\cite{Dohrmann-2003-PSC}, and recently by \v{S}\'{i}stek et al.~\cite{Sistek-2012-FSC}.

Moreover, solving the coarse problem exactly in the original BDDC method
becomes a bottleneck as the number of unknowns and, in particular, the
number of substructures gets too large. Since the coarse problem in BDDC,
unlike in the FETI-DP, has the same structure as the original problem, it is
straightforward to apply the method recursively to solve the coarse problem
only approximately \cite{Dohrmann-2003-PSC}. The original, two-level, BDDC\
has been extended into three-levels by Tu \cite{Tu-2007-TBT3D,Tu-2007-TBT}
and into a general multilevel method by Mandel, Soused\'{\i}k and Dohrmann~%
\cite{Mandel-2007-OMB,Mandel-2008-MMB}. Recently the BDDC\ has been extended
into three-level methods for\ mortar discretizations~\cite{Kim-2009-TLB}, and into 
multiple levels for saddle point problems~\cite{Sousedik-2011-NBS,Tu-2011-TBA}. 
The abstract condition number
bounds deteriorate exponentially with increasing number of levels. 

Here we combine the adaptive and multilevel approaches to the BDDC method in
order to develop its variant that would preserve parallel scalability with
an increasing number of subdomains and also show its excellent convergence
properties. The adaptive method works as previously. It selects constraints
associated with substructure faces, obtained from solution of local
generalized eigenvalue problems for pairs of adjacent substructures, however
this time on each decomposition level. Because of the multilevel approach,
the coarse problems are treated explicitly (unlike in \cite%
{Mandel-2012-ABT,Sousedik-2008-CDD}). The numerical examples show that the
heuristic eigenvalue-based estimates work reasonably well and that the
adaptive approach can result on each decomposition level in the
concentration of computational work in a small troublesome parts of the
problem, which leads to a good convergence behavior.
The developed open-source parallel implementation shows a good scalability as well as 
applicability to very large problems and core counts.

The theoretical part of this paper presents a part of the work from the thesis~\cite{Sousedik-2010-AMB-thesis} in a shorter, 
self-contained way.
Also some results by the serial implementation of the algorithm from \cite{Sousedik-2010-AMB-thesis} are reproduced here  
for comparisons.
The two-dimensional version of the algorithm was described in conference proceedings~\cite{Sousedik-2011-AMB}.
The main original contribution of this paper is the description of the parallel implementation of the method,
and numerical study of its performance.

The paper is organized as follows. In Section~\ref{sec:FE-setting} we
establish the notation and introduce problem settings and preliminaries. In
Section~\ref{sec:multilevel-bddc} we recall the Multilevel BDDC originally
introduced in \cite{Mandel-2008-MMB}. In Section \ref{sec:adaptive}, we
describe the adaptive two-level method in terms of the BDDC\ operators with
an explicit coarse space. In Section~\ref{sec:lobpcg} we discuss a\
preconditioner for LOBPCG used in the solution of the local generalized
eigenvalue problems in the adaptive method. Section~\ref{sec:AM-BDDC}
contains an algorithm for the adaptive selection of components of the
Multilevel BDDC preconditioner. Numerical results are presented in Section~%
\ref{sec:numerical}, and Section~\ref{sec:conclusion} contains summary and
concluding remarks.

\section{Notation and substructuring components}

\label{sec:FE-setting}

We first establish notation and briefly review standard substructuring
concepts and describe BDDC components.
See, e.g.,~\cite{Kruis-2006-DDM,Smith-1996-DD,Toselli-2005-DDM} 
for more details about iterative substructuring in general, and in particular 
\cite{Dohrmann-2003-PSC,Mandel-2003-CBD,Mandel-2007-BFM,Mandel-2008-MMB} 
for the BDDC. Consider a bounded domain $\Omega \subset \mathbb{R}^{3}$ discretized
by conforming finite elements. The domain $\Omega $ is decomposed into $N$
nonoverlapping \emph{subdomains} $\Omega ^{i}$, $i=1,\dots N$, also called 
\emph{substructures}, so that each substructure $\Omega ^{i}$ is a union of
finite elements. Each node is associated with one degree of freedom in the
scalar case, and with 3 displacement degrees of freedom in the case of
linear elasticity. The nodes contained in the intersection of at least two
substructures are called boundary nodes. The union of all boundary nodes of
all substructures is called the \emph{interface}, denoted by $\Gamma $, and $%
\Gamma ^{i}$ is the interface of substructure~$\Omega ^{i}$. The interface $%
\Gamma $ may also be classified as the union of three different types of
sets: \emph{faces}, \emph{edges} and \emph{corners}. We will adopt here a
simple (geometric) definition: a \emph{face} contains all nodes shared by
the same two subdomains, an \emph{edge} contains all nodes shared by same
set of more than two subdomains, and a \emph{corner} is a degenerate edge
with only one node; for a more general definition see, e.g.,~\cite%
{Klawonn-2006-DPF}. 

We identify finite element functions with the vectors of their coefficients
in the standard finite element basis. These coefficients are also called 
\emph{variables} or \emph{degrees of freedom}. We also identify linear
operators with their matrices, in bases that will be clear from the context.

Here, we find it more convenient to use the notation of abstract linear
spaces and linear operators between them instead of the space $\mathbb{R}%
^{n} $ and matrices. The results can be easily converted to matrix language
by choosing a finite element basis. The space of the finite element
functions on $\Omega $ will be denoted as $U$. Let $W^{s}$ be the space of
finite element functions\ on substructure$~\Omega ^{s}$, such that all of
their degrees of freedom on $\partial \Omega ^{s}\cap \partial \Omega $ are
zero. Let%
\begin{equation*}
W=W^{1}\times \cdots \times W^{N},
\end{equation*}%
and consider a bilinear form $a\left( \cdot ,\cdot \right) $ arising from
the second-order elliptic problem such as Poisson's equation or a problem of
linear elasticity.

Now $U\subset W$ is the subspace of all functions from $W$ that are
continuous across the substructure interfaces. We are interested in the
solution of the problem
\begin{equation}
u\in U:a(u,v)=\left\langle f,v\right\rangle ,\quad\forall v\in U,
\label{eq:problem-full-space}
\end{equation}
where the bilinear form $a$ is associated on the space $U$\ with the system
operator$~A$, defined by 
\begin{equation}
A:U\mapsto U^{\prime},\quad a(u,v)=\left\langle Au,v\right\rangle, \quad \forall
u,v\in U,  \label{eq:def-A}
\end{equation}
and $f\in U^{\prime}$ is the right-hand side. Hence, (\ref%
{eq:problem-full-space}) is equivalent to%
\begin{equation}
Au=f.  \label{eq:problem-full-space-algebraic}
\end{equation}

Define $U_{I}\subset U$ as the subspace of functions that are zero on the
interface~$\Gamma$, 
i.e., the `interior' functions. Denote by 
$P$ the energy orthogonal projection from $W$ onto $U_{I}$,%
\begin{equation*}
P:w\in W\longmapsto v_{I}\in U_{I}:a\left( v_{I},z_{I}\right) =a\left(
w,z_{I}\right) ,\quad\forall z_{I}\in U_{I}.
\end{equation*}
Functions from $\left( I-P\right) W$, i.e., from the nullspace of $P,$ are
called discrete harmonic; these functions are $a$-orthogonal to $U_{I}$ and
energy minimal with respect to increments in $U_{I}$. Next, let $\widehat{W}$
be the space of all discrete harmonic functions that are continuous across
substructure boundaries, that is 
\begin{equation}
\widehat{W}=\left( I-P\right) U.  \label{eq:discrete-harm}
\end{equation}
In particular, 
\begin{equation}
U=U_{I}\oplus\widehat{W},\quad U_{I}\perp_{a}\widehat{W}.
\label{eq:int-harm-dec}
\end{equation}

The BDDC\ method \cite{Dohrmann-2003-PSC,Mandel-2007-BFM} is a two-level
preconditioner characterized by the selection of certain \emph{coarse
degrees of freedom}, such as values at the corners and averages over edges
or faces of substructures. Define $\widetilde{W}\subset W$ as the subspace
of all functions such that the values of any coarse degrees of freedom have
a common value for all relevant substructures and vanish on $\partial \Omega
,$ and $\widetilde{W}_{\Delta }\subset \widetilde{W}$ as the subspace of all
functions such that their coarse degrees of freedom vanish. Next, define $%
\widetilde{W}_{\Pi }$ as the subspace of all functions such that their
coarse degrees of freedom between adjacent substructures coincide, and such
that their energy is minimal. Clearly, functions in $\widetilde{W}_{\Pi }$
are uniquely determined by the values of their coarse degrees of freedom,
and 
\begin{equation}
\widetilde{W}_{\Delta }\perp _{a}\widetilde{W}_{\Pi },\text{\quad and\quad }%
\widetilde{W}=\widetilde{W}_{\Delta }\oplus \widetilde{W}_{\Pi }.
\label{eq:tilde-dec}
\end{equation}%
The component of the BDDC\ preconditioner formulated in the space $%
\widetilde{W}_{\Pi }$ is called the \emph{coarse problem} and the components in the
space $\widetilde{W}_{\Delta }$ are called \emph{substructure corrections}.

We assume that 
\begin{equation}
a\text{ is positive definite on }\widetilde{W}.  \label{eq:pos-def}
\end{equation}%
This will be the case when $a$ is positive definite on the space $U$
and there are sufficiently many
coarse degrees of freedom \cite{Mandel-2008-MMB}. We further assume that the coarse degrees of
freedom are zero on all functions from $U_{I}$, that is,%
\begin{equation}
U_{I}\subset \widetilde{W}_{\Delta }.  \label{eq:coarse-int}
\end{equation}%
In other words, the coarse degrees of freedom depend on the values on
substructure boundaries only. From (\ref{eq:tilde-dec}) and (\ref%
{eq:coarse-int}), it follows that the functions in~$\widetilde{W}_{\Pi }$
are discrete harmonic, that is,%
\begin{equation}
\widetilde{W}_{\Pi }=\left( I-P\right) \widetilde{W}_{\Pi }.
\label{eq:coarse-is-discrete-harmonic}
\end{equation}

Next, let $E$ be a projection from $\widetilde{W}$ onto $U$, defined by
taking some weighted average on substructure interfaces. That is, we assume
that%
\begin{equation}
E:\widetilde{W}\rightarrow U,\quad EU=U,\quad E^{2}=E.  \label{eq:E-onto-U}
\end{equation}%
Since a projection is the identity on its range, it follows that $E$ does
not change the interior degrees of freedom, 
\begin{equation}
EU_{I}=U_{I},  \label{eq:int-unchanged}
\end{equation}%
since $U_{I}\subset U$. Finally, we recall that the operator $\left(
I-P\right) E$ is a projection~\cite{Mandel-2008-MMB}.

\section{Multilevel BDDC}

\label{sec:multilevel-bddc}

We recall Multilevel BDDC\ which has been introduced as a particular
instance of Multispace BDDC\ in \cite{Mandel-2008-MMB}. The substructuring
components from Section~\ref{sec:FE-setting} will be denoted by an
additional subscript~$_{1},$ as~$\Omega _{1}^{s},$ $s=1,\ldots N_{1}$, etc.,
and called level~$1$. The level~$1$ coarse problem will be called the level~$%
2$ problem. It has the same finite element structure as the original problem
(\ref{eq:problem-full-space}) on level~$1$, so we put $U_{2}=\widetilde{W}%
_{\Pi 1}$. Level~$1$ substructures are level~$2$ elements and level~$1$
coarse degrees of freedom are level~$2$ degrees of freedom. Repeating this
process recursively, level~$i-1$ substructures become level~$i$ elements,
and the level~$i$ substructures are agglomerates of level~$i$ elements.
Level~$i$ substructures are denoted by $\Omega _{i}^{s},$ $s=1,\ldots
,N_{i}, $ and they are assumed to form a conforming triangulation with a
characteristic substructure size $H_{i}$. For convenience, we denote by $%
\Omega _{0}^{s}$ the original finite elements and put $H_{0}=h$. The
interface$~\Gamma _{i}$\ on level$~i$ is defined as the union of all level$%
~i $ boundary nodes, i.e., nodes shared by at least two level$~i$
substructures, and we note that $\Gamma _{i}\subset \Gamma _{i-1}$. Level $%
i-1$ coarse degrees of freedom become level $i$ degrees of freedom. The
shape functions on level~$i$ are determined by minimization of energy with
respect to level~$i-1$ shape functions, subject to the value of exactly one
level $i$ degree of freedom being one and the other level $i$ degrees of
freedom being zero. The minimization is done on each level $i$ element
(level $i-1$ substructure) separately, so the values of level $i-1$ degrees
of freedom are in general discontinuous between level $i-1$ substructures,
and only the values of level $i$\ degrees of freedom between neighbouring
level~$i$\ elements coincide. 

The development of the spaces on level $i$ now parallels the finite element
setting in Section \ref{sec:FE-setting}. Denote $U_{i}=\widetilde{W}_{\Pi,
i-1}$. Let $W_{i}^{s}$ be the space of functions\ on the substructure $%
\Omega_{i}^{s}$, such that all of their degrees of freedom on $\partial
\Omega_{i}^{s}\cap\partial\Omega$ are zero, and let%
\begin{equation*}
W_{i}=W_{i}^{1}\times\cdots\times W_{i}^{N_{i}}.
\end{equation*}
Then $U_{i}\subset W_{i}$ is the subspace of all functions from $W_{i}$ that are
continuous across the interfaces $\Gamma_{i}$. Define $U_{Ii}\subset U_{i}$
as the subspace of functions that are zero on$~\Gamma_{i}$, i.e., the
functions `interior' to the level$~i$
substructures. Denote by $P_{i}$ the energy orthogonal projection from $%
W_{i} $ onto $U_{Ii}$,%
\begin{equation*}
P_{i}:w_{i}\in W_{i}\longmapsto v_{Ii}\in U_{Ii}:a\left(
v_{Ii},z_{Ii}\right) =a\left( w_{i},z_{Ii}\right) ,\quad\forall z_{Ii}\in
U_{Ii}.
\end{equation*}
Functions from $\left( I-P_{i}\right) W_{i}$, i.e., from the nullspace of $%
P_{i},$ are called discrete harmonic on level $i$; these functions are $a$%
-orthogonal to $U_{Ii}$ and energy minimal with respect to increments in $%
U_{Ii}$. Denote by $\widehat{W}_{i}\subset U_{i}$ the subspace of discrete
harmonic functions on level$~i$, that is 
\begin{equation}
\widehat{W}_{i}=\left( I-P_{i}\right) U_{i}.  \label{eq:discrete-harm-ML}
\end{equation}
In particular, $U_{Ii}\perp_{a}\widehat{W}_{i}$. Define $\widetilde{W}%
_{i}\subset W_{i}$ as the subspace of all functions 
such that each coarse degree of freedom on level$~i$ has a common value for all relevant level$~i$\ substructures, 
and $\widetilde{W}_{\Delta i}\subset \widetilde{W}_{i}$ as the subspace
of all functions such that their level $i$ coarse degrees of freedom have zero value.
Define $\widetilde{W}_{\Pi i}$ as the subspace of all functions such that
their level $i$\ coarse degrees of freedom between adjacent substructures
coincide, and such that their energy is minimal. Clearly, functions in $%
\widetilde{W}_{\Pi i}$ are uniquely determined by the values of their level $%
i$\ coarse degrees of freedom, and 
\begin{equation}
\widetilde{W}_{\Delta i}\perp_{a}\widetilde{W}_{\Pi i},\text{\quad}%
\widetilde{W}_{i}=\widetilde{W}_{\Delta i}\oplus\widetilde{W}_{\Pi i}.
\label{eq:tilde-dec-ML}
\end{equation}
We assume that the level$~i$\ coarse degrees of freedom are zero on all
functions from $U_{Ii}$, that is,%
\begin{equation}
U_{Ii}\subset\widetilde{W}_{\Delta i}.  \label{eq:coarse-int-ML}
\end{equation}
In other words, level $i$ coarse degrees of freedom depend on the values on
level$~i$ substructure boundaries only. From (\ref{eq:tilde-dec-ML}) and (%
\ref{eq:coarse-int-ML}), it follows that the functions in $\widetilde{W}%
_{\Pi i}$ are discrete harmonic on level$~i$, that is%
\begin{equation}
\widetilde{W}_{\Pi i}=\left( I-P_{i}\right) \widetilde{W}_{\Pi i}.
\label{eq:coarse-is-discrete-harmonic-ML}
\end{equation}
Let $E$ be a projection from $\widetilde{W}_{i}$ onto $U_{i}$, defined by
taking some weighted average on $\Gamma_{i}$%
\begin{equation*}
E_{i}:\widetilde{W}_{i}\rightarrow U_{i},\quad
E_{i}^{2}=E_{i}.
\end{equation*}
Since projection is the identity on its range,\ $E_{i}$ does not change the
level$~i$ interior degrees of freedom, in particular%
\begin{equation}
E_{i}U_{Ii}=U_{Ii}.  \label{eq:int-unchanged-ML}
\end{equation}

The Multilevel BDDC\ method is now defined recursively \cite%
{Dohrmann-2003-PSC,Mandel-2008-MMB} by solving the coarse problem on level $%
i $ only approximately, by one application of the preconditioner on level $%
i+1$. Eventually, at the top level $L-1$, the coarse problem, which is the
level~$L$ problem, is solved exactly. A formal description of the method is
provided by the following algorithm.

\begin{algorithm}[{Multilevel BDDC, \protect\cite[Algorithm~17]%
{Mandel-2008-MMB}}]
\label{alg:multilevel-bddc}Define the preconditioner $r_{1}\in U_{1}^{\prime
}\longmapsto u_{1}\in U_{1}$ as follows:

\noindent\textbf{for }$i=1,\ldots,L-1$\textbf{,}

\begin{description}
\item Compute interior pre-correction on level $i$,%
\begin{equation}
u_{Ii}\in U_{Ii}:a\left( u_{Ii},z_{Ii}\right) =\left\langle
r_{i},z_{Ii}\right\rangle ,\quad\forall z_{Ii}\in U_{Ii}.  \label{eq:ML-uIi}
\end{equation}

\item Get an updated residual on level $i$,%
\begin{equation}
r_{Bi}\in U_{i},\quad\left\langle r_{Bi},v_{i}\right\rangle =\left\langle
r_{i},v_{i}\right\rangle -a\left( u_{Ii},v_{i}\right) ,\quad\forall v_{i}\in
U_{i}.  \label{eq:ML-rBi}
\end{equation}

\item Find the substructure correction on level $i$: 
\begin{equation}
w_{\Delta i}\in W_{\Delta i}:a\left( w_{\Delta i},z_{\Delta i}\right)
=\left\langle r_{Bi},E_{i}z_{\Delta i}\right\rangle ,\quad\forall z_{\Delta
i}\in W_{\Delta i}.  \label{eq:ML-wDi}
\end{equation}

\item Formulate the coarse problem on level $i$, 
\begin{equation}
w_{\Pi i}\in W_{\Pi i}:a\left( w_{\Pi i},z_{\Pi i}\right) =\left\langle
r_{Bi},E_{i}z_{\Pi i}\right\rangle ,\quad\forall z_{\Pi i}\in W_{\Pi i}.
\label{eq:ML-coarse}
\end{equation}

\item If\ $\ i=L-1$, solve the coarse problem directly and set $u_{L}=w_{\Pi
L-1}$, \newline
otherwise set up the right-hand side for level $i+1$,%
\begin{equation}
r_{i+1}\in\widetilde{W}_{\Pi i}^{\prime},\quad\left\langle
r_{i+1},z_{i+1}\right\rangle =\left\langle r_{Bi},E_{i}z_{i+1}\right\rangle
,\quad\forall z_{i+1}\in\widetilde{W}_{\Pi i}=U_{i+1},  \label{eq:ML-ri+1}
\end{equation}
\end{description}

\noindent\textbf{end.}

\medskip

\noindent\textbf{for }$i=L-1,\ldots,1\mathbf{,}$\textbf{\ }

\begin{description}
\item Average the approximate corrections on substructure interfaces on
level $i$,%
\begin{equation}
u_{Bi}=E_{i}\left( w_{\Delta i}+u_{i+1}\right) .  \label{eq:ML-uBi-1}
\end{equation}

\item Compute the interior post-correction on level $i$,%
\begin{equation}
v_{Ii}\in U_{Ii}:a\left( v_{Ii},z_{Ii}\right) =a\left( u_{Bi},z_{Ii}\right)
,\quad\forall z_{Ii}\in U_{Ii}.  \label{eq:ML-vIi}
\end{equation}

\item Apply the combined corrections, 
\begin{equation}
u_{i}=u_{Ii}+u_{Bi}-v_{Ii}.  \label{eq:ML-ui}
\end{equation}
\end{description}

\noindent \textbf{end.}
\end{algorithm}

The condition number bound for Multilevel BDDC is given as follows. 

\begin{lemma}[{\protect\cite[Lemma 20]{Mandel-2008-MMB}}]
\label{thm:def-omega-mlevel}The condition number $\kappa $
of Multilevel BDDC from Algorithm~\ref{alg:multilevel-bddc} satisfies%
\begin{equation}
\kappa \le \omega \equiv \Pi _{i=1}^{L-1}\omega _{i}\,\,,\qquad   \label{eq:def-omega-mlevel}
\end{equation}%
where%
\begin{equation}
\label{eq:def-omega}
\omega _{i}=\sup_{w\in \widetilde{W}_{i}}\frac{
\left\Vert \left(I-P_{i}\right) E_{i}w\right\Vert _{a}^{2}
}{
\left\Vert w\right\Vert
_{a}^{2}}.
\end{equation}
\end{lemma}

For the purpose of the adaptive selection of constraints, 
we use the bound based on jump at the interface defined 
on the subspace of discrete harmonic functions from $\widetilde{W}_{i}$.
More precisely, we modify (\ref{eq:def-omega}) 
using the identity 
\begin{equation}
\label{eq:projection identity}
(I-P_{i})E_{i}(I-P_{i}) = (I-P_{i})E_{i}
\end{equation}
and the fact that $P_{i}$ is an $a$-orthogonal projection, as
\begin{equation}
\label{eq:def-omega-mod}
\begin{aligned}
\omega _{i}&=\sup_{w\in \widetilde{W}_{i}}\frac{
\left\Vert \left(I-P_{i}\right) E_{i}w\right\Vert _{a}^{2}
}{
\left\Vert w\right\Vert
_{a}^{2}} 
=\sup_{w\in \widetilde{W}_{i}}\frac{
\left\Vert 
\left(I-P_{i}\right) E_{i}\left(I-P_{i}\right)w\right\Vert _{a}^{2}
}{
\left\Vert w\right\Vert
_{a}^{2}} \\
&=\sup_{w\in \widetilde{W}_{i}}\frac{
\left\Vert 
\left(I-P_{i}\right) E_{i}\left(I-P_{i}\right)w\right\Vert _{a}^{2}
}{
\left\Vert P_{i}w\right\Vert_{a}^{2}
+
\left\Vert \left(I-P_{i}\right)w\right\Vert_{a}^{2}
}
=\sup_{w\in \widetilde{W}_{i}}\frac{
\left\Vert 
\left(I-P_{i}\right) E_{i}\left(I-P_{i}\right)w\right\Vert _{a}^{2}
}{
\left\Vert \left(I-P_{i}\right)w\right\Vert_{a}^{2}
} \\
&=\sup_{w\in \left(I-P_{i}\right)\widetilde{W}_{i}}\frac{
\left\Vert 
\left(I-P_{i}\right) E_{i}w\right\Vert _{a}^{2}
}{
\left\Vert w\right\Vert_{a}^{2}
}
=\sup_{w\in \left(I-P_{i}\right)\widetilde{W}_{i}}\frac{
\left\Vert 
\left(I - \left(I-P_{i}\right) E_{i}\right)w\right\Vert _{a}^{2}
}{
\left\Vert w\right\Vert_{a}^{2}
}.
\end{aligned}
\end{equation}

The last equality in~(\ref{eq:def-omega-mod}) holds because $\left(I-P_{i}\right)E_{i}$ is a projection 
and the norm of a nontrivial projection in an inner product space depends only on the angle 
between its range and its nullspace~\cite{Ipsen-1995-ABC}. 

\section{Adaptive Coarse Degrees of Freedom}

\label{sec:adaptive}

To simplify notation, we formulate the algorithm for the adaptive selection of the
coarse degrees of freedom for one level at a time and drop the subscript~$i$.
The basic idea of the method is
still the same as in \cite{Mandel-2007-ASF,Mandel-2012-ABT,Sousedik-2008-CDD}.
However, the current formulation in terms of the BDDC\ method, though equivalent and written similarly as in~\cite{Mandel-2012-ABT},
is different enough to allow for an
explicit treatment of the coarse space correction. Therefore, it is suitable for multilevel
extension which will be introduced later in Section~\ref{sec:AM-BDDC}.

As mentioned before, the space $\widetilde{W}$ is constructed using coarse degrees of freedom.
These can be, e.g., values at corners, and averages over edges or faces. The
space $\widetilde{W}$ is then given by the requirement that the coarse
degrees of freedom on adjacent substructures coincide; for this reason, the
terms coarse degrees of freedom and constraints are used interchangeably.
The edge (or face) averages are necessary in 3D problems to obtain
scalability with subdomain size. Ideally, one can prove the polylogarithmic
condition number bound 
\begin{equation}
\kappa\leq const\left( 1+\log\frac{H}{h}\right) ^{2},  \label{eq:polylog}
\end{equation}
where $H$ is the subdomain size and $h$ is the finite element size.

\begin{remark}
\label{rem:polylog}The initial selection of constraints in the proposed
adaptive approach will be done in a way such that~(\ref{eq:polylog}) is satisfied
for problems with sufficiently regular structure.
See, e.g.,~\cite{Klawonn-2002-DPF} for a theoretical justification. 
\end{remark}

To choose the space $\widetilde{W}$, cf.~\cite[Section 2.3]{Mandel-2007-ASF}%
, suppose we are given a space~$X$ and a linear operator $C:W\rightarrow X$ and define, 
\begin{equation}
\widetilde{W}=\left\{ w\in W:C\left( I-E\right) w=0\right\} .
\label{eq:def-Wtilde-dual}
\end{equation}
The values $Cw$ will be called local coarse degrees of freedom, and the space%
$~\widetilde{W}$ consists of all functions $w$ whose local coarse degrees of
freedom on adjacent substructures have zero jumps. To represent their common
values, i.e., the global coarse degrees of freedom of vectors $u\in%
\widetilde{W}$, we use a space $U_{c}$ and a one-to-one linear operator $%
R_{c}:U_{c}\rightarrow X$ such that
\begin{equation*}
\widetilde{W}=\left\{ w\in W:\exists u_{c}\in U_{c}:Cw=R_{c}u_{c}\right\}.
\end{equation*}

Observe that $\left(
I-E\right) Pv=0$ for all $v\in W$, so we can define the space $\widetilde{W}$
in (\ref{eq:def-Wtilde-dual}) using discrete harmonic functions $w\in W_{\Gamma } = \left(
I-P\right) W$, for which 
\begin{equation}
\left( I-\left( I-P\right) E\right) w=\left( I-P\right) \left( I-E\right) w.
\label{eq:IPE}
\end{equation}%
Let us denote $\widetilde{W}_{\Gamma} = (I-P)\widetilde{W} = \widetilde{W} \cap W_{\Gamma }$.
Then the bound (\ref{eq:def-omega}) in the form of the last term in equation (\ref{eq:def-omega-mod}) 
can be found, 
for a fixed level~$i$, as a maximum eigenvalue of an associated
eigenvalue problem, which can be using (\ref{eq:IPE}) written as 
\begin{equation}
a\left( \left( I-P\right) \left( I-E\right) w,\left( I-P\right) \left(
I-E\right) z\right)=\lambda a\left( w,z\right)
\quad \forall z\in \widetilde{W}_{\Gamma }.  \label{eq:global-eig-var}
\end{equation}

We can then control the condition number bound by adding constraints 
adaptively by taking advantage of the Courant-Fisher-Weyl minimax principle, cf., e.g., 
\cite[%
Theorem 5.2]{Demmel-1997-ANL}, in the same way as in \cite{Mandel-2007-ASF,Mandel-2012-ABT,Sousedik-2008-CDD}.

\begin{corollary}[\protect\cite{Mandel-2012-ABT}]
\label{cor:optimal-eig}The generalized eigenvalue problem (\ref%
{eq:global-eig-var}) has eigenvalues $\lambda _{1}\geq \lambda _{2}\geq
\ldots \geq \lambda _{n}\geq 0$. Denote the corresponding eigenvectors by $%
w_{\ell }$. Then, for any $k=1,\ldots ,n-1$, and any linear functionals $%
L_{\ell }$, 
$\ell =1,\ldots ,k$,%
\begin{equation*}
\max \left\{ \frac{\left\Vert \left( I-P \right) \left( I-E \right)
w\right\Vert _{a}^{2}}{\left\Vert w\right\Vert _{a}^{2}}:w\in \widetilde{W}%
_{\Gamma },\text{ }L_{\ell }\left( w\right) =0\text{ }\forall \ell =1,\ldots
,k\right\} \geq \lambda _{k+1},
\end{equation*}%
with equality if%
\begin{equation}
L_{\ell }\left( w\right) =a\left( \left( I-P\right) \left( I-E\right)
w_{\ell },\left( I-P\right) \left( I-E\right) w\right) .
\label{eq:extra-constr-var}
\end{equation}
\end{corollary}

Therefore, because $\left( I-E\right) $ is a projection, the optimal
decrease of the condition number bound (\ref{eq:def-omega}) can be achieved
by adding to the constraint matrix $C$ in the definition of $\widetilde{W}$
the rows $c_{\ell}$ defined by $c_{\ell}^{T}w=L_{\ell}\left(
w\right) $. 

Solving the global eigenvalue problem (\ref{eq:global-eig-var}) is expensive, 
and the vectors~$c_{\ell}$ are not of the
form required for substructuring, i.e., each $c_{\ell}$ with nonzero entries
corresponding to only one corner, an edge or a face at a time.
For these reasons, we replace (\ref{eq:global-eig-var}) by a collection of
local problems, each defined by considering only two adjacent subdomains 
$\Omega ^{s}$ and $\Omega ^{t}$. 
Subdomains are called adjacent if they share a face. 
All quantities associated with such pairs will be denoted by a superscript $^{st}$. 
In particular, we define
\begin{equation}
W^{st} = W^{s} \times W^{t}, \quad W^{st}_{\Gamma} = (I-P^{st})W^{st},
\end{equation}
where $(I-P^{st})$ realizes the discrete harmonic extension from the local interfaces $\Gamma ^{s}$ and $\Gamma ^{t}$
to interiors.
Thus, functions from $W^{st}_{\Gamma}$ are fully determined by their values at the local interfaces 
$\Gamma ^{s}$ and $\Gamma ^{t}$,
and they may be discontinuous at the common part $\Gamma^{st} = \Gamma^{s} \cap \Gamma^{t}$.

The bilinear form $a^{st}(\cdot,\cdot)$ is associated on the space $W^{st}_{\Gamma}$ with  
the operator~$S^{st}$ of Schur complement with respect to the local interfaces, 
defined by
\begin{equation}
S^{st}:W^{st}_{\Gamma}\mapsto{W^{st}_{\Gamma}}^{\prime},  \quad a^{st}(u,v)=\langle S^{st}u,v\rangle, \quad \forall u,v \in W^{st}_{\Gamma}.  \label{eq:Schur-definition}
\end{equation}
Operator $S^{st}$ is represented by a block-diagonal matrix composed 
of symmetric positive semi-definite matrices $S^{s}$ and $S^{t}$ 
of individual Schur complements of the subdomain matrices with respect to local interfaces $\Gamma^{s}$ and $\Gamma^{t}$,
resp., 
\begin{equation}
\quad S^{st}=\left[ 
\begin{array}{cc}
S^{s} &  \\ 
& S^{t}%
\end{array}
\right].  
\label{eq:global-S}
\end{equation}

The action of the local projection operator $E^{st}$ is realized as a (weighted) average at $\Gamma^{st}$ 
and as an identity operator at $(\Gamma^{s} \cup \Gamma^{t})\backslash \Gamma^{st}$.

Let $C^{st}$ be the operator defining the initial coarse degrees of freedom that are common to both subdomains 
of the pair.
We define the local space of functions with the shared coarse degrees of freedom continuous as
\begin{equation}
\widetilde{W}^{st} = \left\{ w\in W^{st}: \ \ C^{st}(I-E^{st})w = 0 \right\}.
\label{eq:local-tildeWst}
\end{equation}
Finally, we introduce the space $\widetilde{W}^{st}_{\Gamma} = \widetilde{W}^{st} \cap W^{st}_{\Gamma}$.

Now the generalized eigenvalue problem (\ref{eq:global-eig-var})
becomes a \emph{localized} problem to find $w\in \widetilde{W}_{\Gamma}^{st}$ such that 
\begin{equation}
a^{st}\left( \left( I-P^{st}\right) \left( I-E^{st}\right) w,\left(
I-P^{st}\right) \left( I-E^{st}\right) z\right) =\lambda \,a^{st}\left(
w,z\right) \quad \forall z\in \widetilde{W}_{\Gamma }^{st}.
\label{eq:global-eig-var-ij}
\end{equation}%

\begin{assumption}
\label{assum:local-regular}The corner constraints are already sufficient to
prevent relative rigid body motions of any pair of adjacent substructures,
so 
\begin{equation*}
\forall w\in \widetilde{W}^{st}:a^{st}(w,w)=0\Rightarrow \left( I-{E}%
^{st}\right) w=0,
\end{equation*}%
i.e., the corner degrees of freedom are sufficient to constrain the rigid
body modes of the two substructures into a single set of rigid body modes,
which are continuous across the interface~$\Gamma ^{st}$.
\end{assumption}

The maximal eigenvalue $\omega^{st}$ of (\ref{eq:global-eig-var-ij}) is
finite due to Assumption \ref{assum:local-regular}, and we define the
heuristic \emph{condition number indicator}
\begin{equation}
\widetilde{\omega}=\max\left\{ \omega^{st}:\Omega^{s}\text{ and }\Omega ^{t}%
\text{ are adjacent}\right\} .  \label{eq:cond-ind}
\end{equation}

Considering two adjacent subdomains $\Omega^{s}$ and $\Omega^{t}$ only, we
get the added constraints $L_{\ell}\left( w\right) =0$ from (\ref%
{eq:extra-constr-var}) as 
\begin{equation}
a^{st}\left( \left( I-P^{st}\right) \left( I-E^{st}\right) w_{\ell },\left(
I-P^{st}\right) \left( I-E^{st}\right) w\right) =0\quad
\forall\ell=1,\ldots,k,  \label{eq:extra-constr-var-ij}
\end{equation}
where $w_{\ell}$ are the eigenvectors corresponding to the $k$ largest
eigenvalues from~(\ref{eq:global-eig-var-ij}).

Let us denote $D$ the matrix corresponding to $C^{st}(I-E^{st})$.
We define the orthogonal projection onto $\mbox{null}\,D$ by 
\begin{equation*}
\Pi=I-D^{T}\left( DD^{T}\right) ^{-1}D.
\end{equation*}
The generalized eigenvalue problem (\ref{eq:global-eig-var}) now becomes
\begin{equation}
\Pi\left( I-P^{st}\right) ^{T}\left( I-E^{st}\right) ^{T}S^{st}\left(
I-E^{st}\right) \left( I-P^{st}\right) \Pi w=\lambda\Pi S^{st}\Pi w.  \label{eq:eig-matrix}
\end{equation}
Since
\begin{equation}
\mbox{null}\Pi S^{st}\Pi\subset\mbox{null}\Pi\left( I-P^{st}\right) ^{T}\left( I-E^{st}\right) ^{T}S^{st}\left( I-E^{st}\right) \left( I-P^{st}\right) \Pi,
\label{eq:eig-nullspace}
\end{equation}
the eigenvalue problem (\ref{eq:eig-matrix}) reduces in the factorspace
modulo $\mbox{null}\Pi S^{st}\Pi$ to a problem with the operator on the
right-hand side positive definite. In our computations, we have used
the subspace iteration method LOBPCG \cite{Knyazev-2001-TOP} to find the
dominant eigenvalues and their eigenvectors. The LOBPCG iterations then
simply run in the factorspace. 

From (\ref{eq:eig-matrix}), the
constraints to be added are 
\begin{equation*}
L_{\ell}\left( w\right) =w_{\ell}^{T}\Pi\left( I-P^{st}\right) ^{T}\left( I-E^{st}\right) ^{T}S^{st}\left( I-E^{st}\right) \left( I-P^{st}\right) \Pi
w=0.
\end{equation*}
That is, we wish to add to the constraint matrix $C$ the rows
\begin{equation}
c^{st}_{\ell}=w_{\ell}^{T}\Pi\left( I-P^{st}\right) ^{T}\left(
I-E^{st}\right) ^{T}S^{st}\left( I-E^{st}\right) \left( I-P^{st}\right) \Pi.
\label{eq:d-def}
\end{equation}

\begin{proposition}[\protect\cite{Mandel-2012-ABT}]
\label{prop:entry-compatible}The vectors $c^{st}_{\ell }$, constructed for a
domain consisting of only two substructures $\Omega ^{s}$ and $\Omega ^{t}$,
have matching entries on the interface between the two substructures, with
opposite signs.
\end{proposition}

That is, each row $c^{st}_{\ell}$ can be split into two blocks and written as
\begin{equation*}
c^{st}_{\ell}=\left[ 
\begin{array}{cc}
c_{\ell}^{s} & \,\,\, -c_{\ell}^{s}%
\end{array}
\right] .
\end{equation*}
Either half of each row from the block $c^{st}_\ell$ is then added into the matrices 
$C^{s}$ and $C^{t}$ corresponding to 
the subdomains $\Omega^{s}$ and $\Omega^{t}$. 
Unfortunately, the added rows will generally have nonzero entries over
the whole $\Gamma^{s}$ and $\Gamma^{t}$, including the edges
in 3D where $\Omega^{s}$ and $\Omega^{t}$ intersect other substructures. 
Consequently,  the added rows are not of the form required for substructuring, i.e., each row 
with nonzeros in one edge or face only. 
In the computations reported in Section~\ref{sec:numerical}, we drop the
adaptively generated edge constraints in 3D. Then it is no longer guaranteed
that the condition number indicator $\widetilde{\omega}\leq\tau$. However,
the method is still observed to perform well.

The proposed adaptive algorithm follows.

\begin{algorithm}[Adaptive BDDC \protect\cite{Mandel-2007-ASF}]
Find the smallest $k$ for every two adjacent substructures $\Omega^{s}$ and $\Omega^{t}$ to guarantee that 
$\lambda^{st}_{k+1}\leq\tau$, where $\tau$ is a given tolerance, and add the
constraints (\ref{eq:extra-constr-var-ij}) to the definition of $\widetilde{W}$.
\end{algorithm}

\section{Preconditioned LOBPCG}

\label{sec:lobpcg}

As pointed out already for \emph{adaptive 2-level BDDC} method in \cite{Mandel-2012-ABT}, 
an important step for a~parallel implementation of the adaptive selection of constraints 
is an efficient solution of the generalized eigenvalue problem (\ref{eq:eig-matrix}) for each pair of adjacent subdomains. 

There are several aspects of the method immediately making such implementation challenging:
(i) parallel layout of pairs of subdomains does not follow the natural layout of a domain decomposition computation with distribution of data based on subdomains,
(ii) the multiplication by $S^{st}$ on both sides of equation (\ref{eq:eig-matrix}) is done only implicitly, 
since action of $S^s$ and $S^t$ is available only through solution of local discrete Dirichlet problems on subdomains $\Omega^s_i$ and~$\Omega^t_i$,
(iii) the process responsible for solving an $st$-eigenproblem typically does not have data for subdomains~$\Omega^s_i$ and~$\Omega^t_i$,
and thus it has to communicate the vector for multiplication to processors able to compute the actions of~$S^s$ and~$S^t$.

With respect to these issues, it is necessary to use an inverse-free method for the solution of each of these problems. 
In our case, the LOBPCG method~\cite{Knyazev-2001-TOP} is applied to find several largest eigenvalues~$\lambda_\ell$ and corresponding eigenvectors~$w_\ell$ solving the homogeneous problem
\begin{equation}
{\mathcal M}({\mathcal A} - \lambda_\ell {\mathcal B})w_\ell = 0,
\label{eq:gen_eig_problem}
\end{equation}
with 
\begin{equation*}
{\mathcal A} = \Pi\left( I-P^{st}\right) ^{T}\left( I-E^{st}\right) ^{T}S^{st}\left( I-E^{st}\right) \left( I-P^{st}\right) \Pi,\quad{\mathcal B} = \Pi S^{st}\Pi,
\end{equation*}
and ${\mathcal M}$ a~suitable preconditioner.
The LOBPCG method requires only multiplications by matrices ${\mathcal M}$, ${\mathcal A}$, and ${\mathcal B}$,
and it can run in the factorspace with~${\mathcal B}$ only positive semi-definite. 
This is important for our situation -- 
although each pair of subdomains has enough initial constraints by corners and edge averages to avoid mechanisms between the two substructures 
(enforced by the projection $\Pi$), 
no essential boundary conditions are applied to the pair as a~whole, and matrix~${\mathcal B}$ of such `floating' pair typically has nontrivial nullspace 
(e.g. rigid body modes for elasticity problems). 

Initial experiments in \cite{Mandel-2012-ABT,Sousedik-2008-CDD,Sousedik-2010-AMB-thesis} 
revealed that while the unpreconditioned LOBPCG (${\mathcal M} = I$) works reasonably
well for simple problems, it requires prohibitively many iterations for problems with very
irregular substructures and/or high jumps in coefficients. 
Since each iteration requires communicating the vector for multiplication,
reducing the iteration counts of LOBPCG by preconditioning is a~very sensible way of accelerating the adaptive BDDC method.

Recall, that the BDDC method provides a preconditioner for the interface problem of the Schur complement by the exact solution of the problem
at the larger space $\widetilde{W}$. 
As such, components of a BDDC implementation, the coarse solver and subdomain corrections, 
can be used to determine the approximate action of the Moore-Penrose pseudoinverse of the matrix ${\mathcal B}$, 
denoted as $M^{loc}_{BDDC}\approx(\Pi S^{st}\Pi)^{+}$. 
This operator can be used as the preconditioner ${\mathcal M}$ for problem (\ref{eq:gen_eig_problem}),
effectively converting the generalized eigenvalue problem to an ordinary eigenproblem inside the iterations.
Using notation from (\ref{eq:global-S}), the preconditioner is formally written as
\begin{equation}
\label{eq:bddc_pair}
M^{loc}_{BDDC} = 
\Pi     
\left(
\left[
\begin{matrix}
I & 0 
\end{matrix}
\right]
\left[
\begin{matrix}
S^{st} & C^T \\
C      & 0  
\end{matrix}
\right]^{-1}
\left[
\begin{matrix}
I \\
0 
\end{matrix}
\right]
+
\Psi ( \Psi ^T S^{st} \Psi )^{+} \Psi ^T
\right)
\Pi,
\end{equation}
where in addition $C = \left[ \begin{matrix} C^s &  \\ & C^t \end{matrix}\right]$ is the matrix of initial constraints 
(continuity at corners and arithmetic averages on edges), 
and $\Psi = \left[ \begin{matrix} \Psi^s R^s_c \\ \Psi^t R^t_c \end{matrix}\right]$ denotes the matrix of coarse basis functions for the two subdomains.
Here $R^i_c,\ i = s,t$, is the zero-one matrix of restriction of the vector of global coarse degrees of freedom of the pair to subdomain coarse degrees of freedom.
Since some coarse degrees of freedom are shared by the two subdomains, corresponding columns in $\Psi$ are nonzero in both parts, 
while columns of coarse degrees of freedom not common to the two subdomains are only nonzero in either $\Psi^s R^s_c$ or $\Psi^t R^t_c $.
Let us recall that in BDDC, the local coarse basis functions are computed as the solution to the problem with multiple right hand sides
\begin{equation}
\label{eq:local_saddle_point_problem}
\left[
\begin{matrix}
S^s &     & C^{sT} &        \\
    & S^t &        & C^{tT} \\
C^s &     &        &        \\
    & C^t &        & 
\end{matrix}
\right]
\left[
\begin{matrix}
\Psi^s &        \\
       & \Psi^t \\
\mu ^s &        \\
       &  \mu ^t
\end{matrix}
\right]
= 
\left[
\begin{matrix}
    &   \\
    &   \\
I^s &   \\
    & I^t
\end{matrix}
\right],
\end{equation}
which represents an independent saddle-point problem with invertible matrix for each subdomain, 
and factorization of which is later reused in applications of the preconditioner (\ref{eq:bddc_pair}).
We also note that the coarse matrix is in the implementation explicitly computed using the second part of the solution 
of (\ref{eq:local_saddle_point_problem}) as 
\begin{equation}
\Psi ^T S^{st} \Psi = - R^{sT}_c \mu ^s R^{s}_c - R^{tT}_c \mu ^t R^{t}_c.
\end{equation}

The coarse matrix of the $st$-pair $\Psi ^T S^{st} \Psi$ has dimension of the union of the coarse degrees of freedom of the subdomains of the pair,
and it is typically only positive semi-definite for a~floating pair. 
Due to its small dimension, we compute its pseudoinverse by means of dense eigenvalue decomposition performed by the LAPACK library,
\begin{equation}
\Psi ^T S^{st} \Psi = V \Lambda V^T, \quad (\Psi ^T S^{st} \Psi )^{+} \approx V \Lambda'^{-1} V^T, 
\end{equation}
where diagonal matrix $\Lambda '$ arises from $\Lambda$ by dropping eigenvalues lower than a prescribed tolerance.

Unlike in the standard BDDC preconditioner, 
no interface averaging is applied to the function before and after the action of $M^{loc}_{BDDC}$,
because problem (\ref{eq:eig-matrix}) is defined in the space $\widetilde{W}^{st}_{\Gamma}$.
Correspondingly, the only approximation is due to using $\Lambda '$ instead of $\Lambda$.

\section{Adaptive-Multilevel BDDC}

\label{sec:AM-BDDC}

We build on Sections~\ref{sec:multilevel-bddc} and \ref{sec:adaptive} to propose a new variant of the
Multilevel BDDC with adaptive selection of constraints on each level.

The development of adaptive selection of constraints in Multilevel BDDC\ now
proceeds similarly as in Section \ref{sec:adaptive}. 
We formulate (\ref{eq:def-omega}) as a set of eigenvalue problems for each
decomposition level. On each level we solve for every two adjacent
substructures a generalized eigenvalue problem and we add the constraints to
the definitions of $\widetilde{W}_{i}$.

The heuristic \emph{condition number indicator }is defined as 
\begin{equation}
\widetilde{\omega}=\Pi_{i=1}^{L-1} \widetilde{\omega}_{i},\quad \widetilde{\omega}_{i} = \max\left\{ \omega_{i}^{st}:\Omega_{i}^{s}%
\text{ and }\Omega_{i}^{t}\text{ are adjacent}\right\} .
\label{eq:cond-ind-mlevel}
\end{equation}

We now describe the \emph{Adaptive-Multilevel BDDC} in more detail. The
algorithm consists of two main steps: (i) set-up (including adaptive selection of
constraints), and (ii) loop of the preconditioned conjugate gradients (PCG) with
the Multilevel BDDC from Algorithm~\ref{alg:multilevel-bddc} as a
preconditioner. The set-up can be summarized as follows (cf.~\cite[Algorithm~4%
]{Sousedik-2010-AMB-thesis} for the 2D case):

\begin{algorithm}[Set-up of Adaptive-Multilevel BDDC]
\label{alg:setup_adaptive}
Adding of coarse degrees of freedom 
to guarantee that the condition number indicator $\widetilde{\omega }\leq
\tau ^{L-1}$, for a given target value$~\tau $:

\medskip

\noindent\textbf{for levels} $i=1:L-1 \mathbf{,}$

\begin{description}
\item Create substructures with roughly the same numbers of degrees of
freedom. 
\medskip

\item Find a set of initial constraints (in particular sufficient number of
corners), and 
set up the BDDC structures for the adaptive algorithm (the next loop over
faces).
\end{description}

\textbf{for all faces} $\mathcal{F}_{i}$ \textbf{on level} $i\mathbf{,} $

\begin{description}
\item \hspace{6mm} Compute the largest local eigenvalues and corresponding
eigenvectors, until the first $m^{st}$ is found such that $%
\lambda_{m^{st}}^{st}\leq\tau$.%
\medskip

\item \hspace{6mm} Compute the constraint weights 
and add these rows to 
the subdomain matrices of constraints $C^s$ and $C^t$. 
\end{description}

\textbf{end.} 

\begin{description}
\item Set-up the BDDC structures for level $i$.
\newline
If the prescribed number of levels is reached, solve the problem directly.
\end{description}

\noindent \textbf{end.}
\end{algorithm}

\section{Implementation remarks}

\label{sec:implementation}

Serial implementation has been developed in Matlab in the thesis \cite{Sousedik-2010-AMB-thesis}. 
Parallel results use the open-source package BDDCML\footnote{ 
\url{http://www.math.cas.cz/~sistek/software/bddcml.html}
} (version 2.0).
This solver is written in Fortran~95 programming language and parallelized using MPI library.
Apart of symmetric positive definite problems studied in this paper, the solver also supports symmetric indefinite and 
general non-symmetric linear systems arising from discretizations of PDEs.

The matrices of the averaging operator $E$ were constructed with entries
proportional to the diagonal entries of the substructure matrices before
elimination of interiors, which is also known as the \emph{stiffness scaling} 
\cite{Klawonn-2008-AFA}.

\subsection{Initial constraints}

Following Remark~\ref{rem:polylog}, in order to satisfy the polylogarithmic
condition number bounds, we have used corners and
arithmetic averages over edges as initial constraints. It is essential
(Assumption~\ref{assum:local-regular}) to generate a sufficient number of
initial constraints to prevent rigid body motions between any
pair of adjacent substructures. 
The selection of corners in our parallel implementation follows the recent face-based algorithm from \cite{Sistek-2012-FSC}.

\subsection{Adaptive constraints}
The adaptive algorithm uses matrices and operators that are readily
available in an implementation of the BDDC method with an explicit coarse
space, with one exception: in order to satisfy the local partition of unity,
cf.~\cite[eq.~(9)]{Mandel-2007-BFM},
\begin{equation*}
E_{i}^{st}R_{i}^{st}=I,
\end{equation*}
we need to generate locally the weight matrices $E_{i}^{st}$ to act as 
an identity operator at $(\Gamma^{s} \cup \Gamma^{t})\backslash \Gamma^{st}$ 
(cf. Section~\ref{sec:adaptive}).

In the computations reported in Section~\ref{sec:numerical}, we drop the
adaptively generated edge constraints in 3D. Then, it is no longer
guaranteed that the condition number indicator $\widetilde{\omega }\leq
\tau^{L-1}$. However, the method is still observed to perform well. Since
the constraint weights are thus supported only on faces, and the entries
corresponding to edges are set to be zero, we orthogonalize and normalize
the vectors of constraint weights (by reduced QR decomposition from LAPACK)
to preserve numerical stability.

In our experience, preconditioning of the LOBPCG method as described in Sec.~\ref{sec:lobpcg} led to a considerable reduction of the number of LOBPCG iterations.  
Or in other words, since we usually put a limit of maximum 15 iterations for an eigenproblem, the resulting eigenvectors are much better converged than 
without preconditioning.
In the parallel implementation of the adaptive selection of constraints, 
pairs are assigned to processors independently of assignment of subdomains.
The BDDCML package uses the open-source implementation of LOBPCG method~\cite{Knyazev-2001-TOP} available in the BLOPEX 
package\footnote{
\url{http://code.google.com/p/blopex}
}. 
Details of the parallel implementation of adaptive selection of constrains were described for 2-level BDDC method in detail in \cite{Sistek-2012-SPA}.

\subsection{Multilevel implementation}

The BDDCML library allows assignment of multiple subdomains at each process. 
At each level, subdomains are assigned to available processors, always starting from root.
Distribution of subdomains on the first level is either provided by user's application, or created by the solver using ParMETIS library (version 3.2). 
On higher levels, where the mesh is considerably smaller, METIS (version 4.0) \cite{Karypis-1998-MSP} is internally used by BDDCML to create mesh partitions.
This means, that on higher levels, where number of subdomains is lower than number of processors, cores with higher ranks are left idle by the preconditioner.

For solving local discrete Dirichlet and Neumann problems on each subdomain, BDDCML relies on a sequential instance of direct solver MUMPS~\cite{Amestoy-2000-MPD}.
A parallel instance of MUMPS is also invoked for factorization and repeated solution of the final coarse problem at the top level.
More details on implementation of the (non-adaptive) multilevel approach in BDDCML can be found in \cite{Sistek-2012-PIM}.

\section{Numerical Examples}

\label{sec:numerical}

To study the properties of the \emph{Adaptive-Multilevel BDDC} method numerically,
we have selected four problems of structural analysis -- two artificial benchmark problems and two realistic engineering problems.
Some of the results were obtained by our serial implementation written in Matlab and reported in thesis~\cite{Sousedik-2010-AMB-thesis}.
This implementation is mainly used to study convergence behaviour with respect to prescribed tolerance on the condition number indicator $\tau$.
The other set of results is obtained using our newly developed parallel implementation within the BDDCML library.
Parallel results were obtained on Cray XE6 supercomputer \emph{Hector} at the Edinburgh Parallel Computing Centre. 
In the computations, one step of the \emph{Adaptive-Multilevel BDDC} method is used as the preconditioner in the preconditioned conjugate gradient (PCG) method, which is run
until the (relative) norm of residual decreases below $10^{-8}$ (in Matlab tests) or $10^{-6}$ (in BDDCML runs).

\subsection{Elasticity in a cube without and with jump in material coefficients}

\begin{figure}[tbph]
\begin{center}
\includegraphics[width=0.36\textwidth]{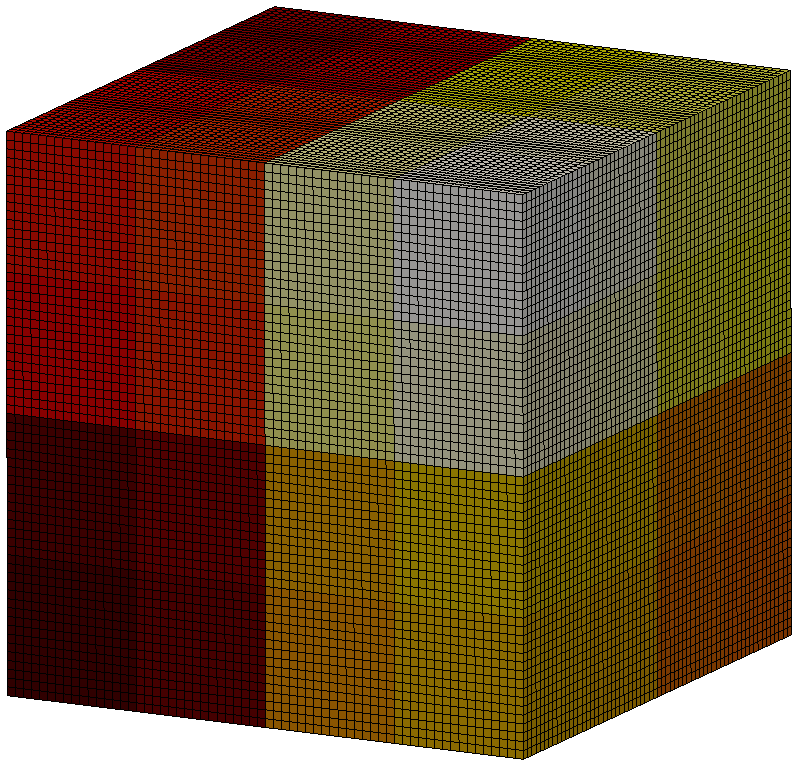}
\hskip 10mm
\includegraphics[width=0.52\textwidth]{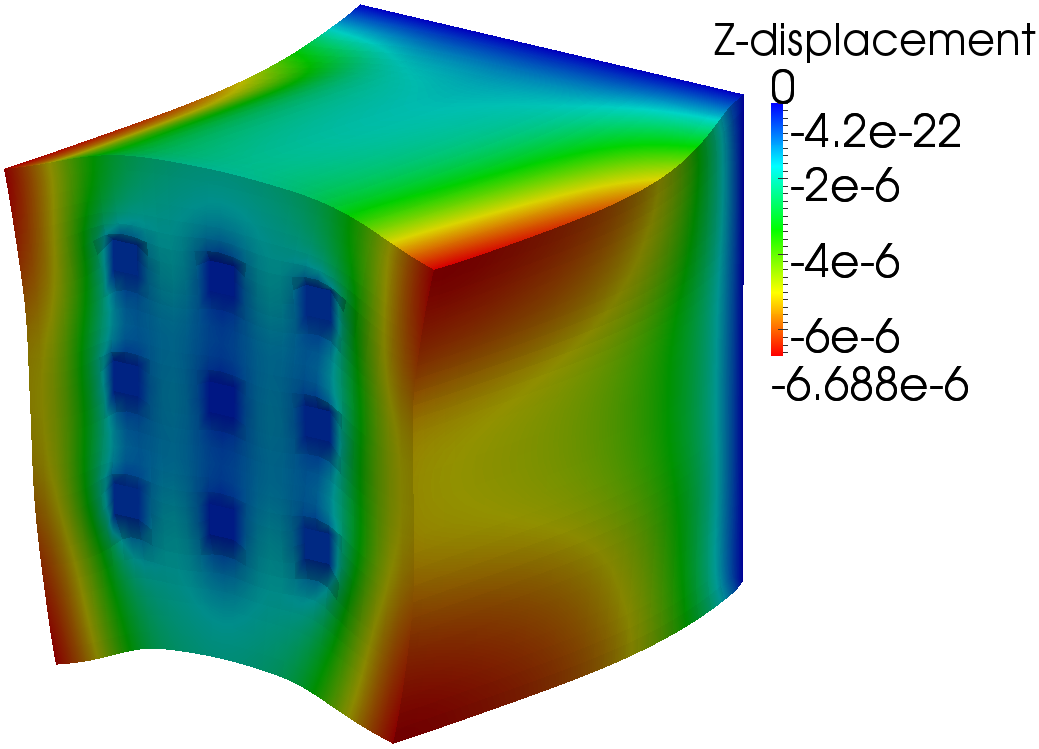}
\caption{\label{fig:cube} Example of a division of the cube into 64 subdomains (left) and (magnified) deformed shape for contrast $E_2/E_1 = 10^6$ coloured by vertical displacement (right).}
\end{center}
\end{figure}

As the first problem, we use the standard benchmark problem of a~unit cube. 
In our setting, we solve the elastic response of the cube under loading by its own weight, when it is fixed at one vertical edge. 
There are nine bars cutting horizontally through the cube. 
We test the case when the bars are of the same material as the rest of the cube (homogeneous material) and
the case when Young's modulus of the outer material $E_1$ is $10^6$ times smaller than that of the bars $E_2$, 
creating contrast in coefficients $E_2/E_1 = 10^6$.
In Fig.~\ref{fig:cube} (right), the (magnified) deformed shape of the cube is shown for this jump in Young's modulus.
We have recently presented a~detailed study of behaviour of the standard (2-level) BDDC method and its adaptive extension with respect to contrast 
on the same problem in \cite{Sistek-2012-SPA}.
It was shown in that reference, that while convergence of BDDC with the standard choice of arithmetic averages on faces quickly deteriorates with increasing contrast,
adaptive version of the algorithm is capable of maintaining good convergence also for large values of contrast, at the cost of quite expensive set-up phase. 

The multilevel approach (without adaptivity), although it may lead to faster solution, suffers from exponentially growing condition number and related number
of iterations,
as reported in \cite{Mandel-2008-MMB}, or recently in \cite{Sistek-2012-PIM}.
Here, we investigate the effect of constraints adaptively generated at higher levels in the multilevel algorithm.
We also study the parallel performance of our solver on this test problem.

The cube is discretized using uniform mesh of tri-linear finite elements and divided into an increasing number of subdomains. 
On the first level, subdomains are cubic with constant $H/h = 16$ ratio 
(see Fig.~\ref{fig:cube} left for an example of a division into 64 subdomains).
On higher levels, divisions into subdomains are created automatically inside BDDCML by the METIS package, 
in general not preserving cubic shape of subdomains.

In Tabs.~\ref{tab:cube_scaling_c1_nonadapt} and \ref{tab:cube_scaling_c1_adapt}, 
we present results of a weak scaling test for the case of the homogeneous cube, i.e. $E_2/E_1 = 1$.
This problem is very well suited for the BDDC method, and the performance is generally very good. 
The growing problem is solved on 8 to 32768 processors (with each core handling one subdomain of the first level).
In these tables, 
$N$ denotes the number of subdomains (and computer cores), 
$n$ denotes global problem size, 
$n_{\Gamma}$ represents the size of the reduced problem defined at the interface $\Gamma$, 
$n_f$ is the number of faces in divisions on the levels (corresponding to number of generalized eigenproblems solved in the adaptive approach), 
`its.' is the number of iterations needed by the PCG method, 
and `cond.' is the estimated condition number obtained from the tridiagonal matrix generated in PCG.
We report times needed by the set-up phase (`set-up'), by PCG iterations (`PCG') and their sum (`solve').

In Tab.~\ref{tab:cube_scaling_c1_nonadapt}, no adaptivity is used, and only the number of levels is varying. 
We can see, that for the standard (2-level) BDDC, we obtain the well-known independence of number of iterations on problem size. 
We can also see, how condition number (and number of PCG iterations) grows when using more levels. 
Although this can lead to savings in time in certain circumstances (due-to cheaper set-up), no such benefits are seen here and these are more common 
in tests of strong scaling with fixed problem size~\cite{Sistek-2012-PIM}.

The independence on problem size is slightly biased on higher levels, probably due to the irregular subdomains.
Computational times slightly grow with problem size, suggesting sub-optimal scaling of BDDCML, especially when going from 512 to 4096 computing cores.
For the largest problem of 32$\times$32$\times$32 subdomains with 405 million degrees of freedom solved on 32768 cores, all times grow considerably. 
This is most likely due to the higher cost of global communication functions at this core count, 
and these results will serve for future performance analysis and optimization of the BDDCML solver. 
Note, that in the case of two levels, parallel direct solver MUMPS failed to solve the resulting coarse problem at this level of parallelism,
which is marked by `n/a' in the tables.

We are now interested in the effect of adaptively generated constraints on convergence of the multilevel BDDC method.
Based on recommendations from \cite{Sistek-2012-SPA}, we limit number of LOBPCG iterations to 15 and maximal number of computed eigenvectors to 10 to 
maintain the cost of LOBPCG solution low.
The target condition number limit is set low, $\tau = 1.5$, which leads to using most of the adaptively generated constraints in actual computation. 
Results are reported in Tab.~\ref{tab:cube_scaling_c1_adapt}. 
We can see, that the adaptive approach is capable of keeping the iteration counts lower, and although the independence of the number of levels is not achieved,
the growth is slower than for the non-adaptive approach.
While the scalability of the solver is similar to the non-adaptive case, 
it is not surprising that the computational time is now dominated by the solution of 
the generalized eigenvalue problems.
This fact makes the adaptive method unsuitable for simple problems like this one, 
in agreement with conclusions for the 2-level BDDC in \cite{Sistek-2012-SPA}.

The situation changes however, when some numerical difficulties appear in the problem of interest. 
One source of such difficulties may be presented by jumps in material coefficients. 
To model this effect, we increased the jump between Young's moduli of the stiff rods and soft outer material to $E_2/E_1 = 10^6$, 
and these results are reported in Tabs.~\ref{tab:cube_scaling_c1e6_nonadapt} and \ref{tab:cube_scaling_c1e6_adapt}. 
For the non-adaptive method (Tab.~\ref{tab:cube_scaling_c1e6_nonadapt}) we can see growth of number of iterations and condition number not only with
adding levels, but also for growing problem size.
This growth is translated to large time spent in PCG iterations, which now dominate the whole solution.

Results are very different for the adaptive approach in Table~\ref{tab:cube_scaling_c1e6_adapt},
for which the main cost is still presented by the solution of the related eigenproblems (included into time of `set-up').
Since we keep the number of computed eigenvectors constant (ten) for each pair of subdomains, the method is not able to maintain a low condition number 
after all these eigenvectors are used for generating constraints.
However, number of iterations is always significantly lower than in the non-adaptive approach, and the method typically requires about one half of the computational time.
While this is an important saving of computational time,
it is also shown in \cite{Sistek-2012-SPA}, that the adaptive approach can solve even problems with contrasts such high, 
that they are not solvable by the non-adaptive approach with arithmetic averages on all faces and edges.

Finally, we compare properties of the coarse basis functions on the first and the second level on this problem.
We consider homogeneous material of the cube which is divided into regular cubic subdomains both on the 
first and the second (unlike in the previous test) level.
Namely, the cube is divided into
4$\times$4$\times$4$=$64 subdomains on the second level.
Each of these subdomains is composed again of 4$\times$4$\times$4$=$64 subdomains of the first level,
which gives 4096 subdomains.
Each of these first-level subdomains is composed of 4$\times$4$\times$4$=$64 tri-linear finite elements.
The problem has in total 262144 elements and 823872 unknowns.
Table~\ref{tab:adaptive_constraints_on_cube} summarizes results of the adaptive 3-level BDDC method 
for different values of prescribed tolerance $\tau$.
For comparison, the non-adaptive 3-level BDDC method with three arithmetic averages on each face requires 19 PCG iterations and the resulting estimated 
condition number is 6.88.

We can see, that significantly (roughly five times) more constraints are selected on the second level than on the first one, which 
suggests that the discrete harmonic basis functions of the first level lead to worse conditioned coarse problem
on the second level.
Thus, it underlines the importance of adaptive selection of constraints on higher levels.
For $\tau ^2 = 2.25$, the maximal number of adaptive constraints (ten) is used on each pair, and
the algorithm is `saturated'. 
Consequently, more constraints would be necessary on each pair to satisfy the condition $\widetilde{\omega}\leq \tau ^2$
from Algorithm~\ref{alg:setup_adaptive}.

\begin{table}[tbph]
\begin{center}
\begin{tabular}
[c]{|c|ccc|cc|ccc|}
\hline
$N$ & 
\multirow{2}{*}{$n$} &
\multirow{2}{*}{$n_{\Gamma}$} &
$n_{f}$ &
\multirow{2}{*}{its.} & 
\multirow{2}{*}{cond.} & 
\multicolumn{3}{c|}{time (sec)} \\
$\ell=1(/2/3)$ &        
&
&        
$\ell=1(/2/3)$ & 
&    
&      
set-up & 
PCG & 
solve \\
\hline
\hline
\multicolumn{9}{|c|}{\textbf{2 levels}} \\
\hline
  8         &  0.1M &  9.5k  &  12            &   15  &  6.7  & 3.9    &  1.6 &    5.5 \\ 
  64        &  0.8M &  0.1M  &  0.1k          &   19  &  7.3  & 4.6    &  2.1 &    6.7 \\
 512        &  6.4M &  1.0M  &  1.3k          &   20  &  6.8  & 9.4    &  3.2 &   12.6 \\
4096        & 50.9M &  8.4M  & 11.5k          &  n/a  &  n/a  & n/a    &  n/a &    n/a \\
\hline
\hline
\multicolumn{9}{|c|}{\textbf{3 levels}} \\
\hline
  64/8      &  0.8M &  0.1M  &  0.1k/18       &   23  &  9.6  &  4.5   &  2.4 &    7.0 \\
 512/64     &  6.4M &  1.0M  &  1.3k/295      &   30  & 16.9  &  5.7   &  3.6 &    9.3 \\
4096/512    & 50.9M &  8.4M  & 11.5k/2930     &   31  & 13.2  & 19.0   &  7.3 &   26.3 \\
32768/128   & 405.0M & 69.1M & 95.2k/664      &   36  & 24.7  & 165.8  & 20.0 &  185.7 \\
\hline
\hline
\multicolumn{9}{|c|}{\textbf{4 levels}} \\
\hline
 512/64/8   &  6.4M &  1.0M  &  1.3k/295/23   &   41  & 24.5  & 5.5    &  4.8 &   10.4 \\
4096/512/64 & 50.9M &  8.4M  & 11.5k/2930/380 &   64  & 87.7  & 9.2    & 11.5 &   20.8 \\
32768/512/8 & 405.0M & 69.1M & 95.2k/2921/23  &   45  & 33.0  & 156.5  & 24.7 &  181.2 \\ 
\hline
\end{tabular}
\end{center}
\caption{\label{tab:cube_scaling_c1_nonadapt} 
Weak scaling for the cube problem with homogeneous material, \emph{non-adaptive multilevel BDDC}.
}
\end{table}

\begin{table}[tbph]
\begin{center}
\begin{tabular}
[c]{|c|ccc|cc|ccc|}
\hline
$N$ & 
\multirow{2}{*}{$n$} &
\multirow{2}{*}{$n_{\Gamma}$} &
$n_{f}$ &
\multirow{2}{*}{its.} & 
\multirow{2}{*}{cond.} & 
\multicolumn{3}{c|}{time (sec)} \\
$\ell=1(/2/3)$ &        
&
&        
$\ell=1(/2/3)$ & 
&    
&      
set-up & 
PCG & 
solve \\
\hline
\hline
\multicolumn{9}{|c|}{\textbf{2 levels}} \\
\hline
  8         &  0.1M &   9.5k &  12            &   11  &  2.5  &   56.1 &  1.2 &   57.3 \\ 
  64        &  0.8M &   0.1M &  0.1k          &   13  &  3.1  &  119.3 &  1.5 &  120.9 \\ 
 512        &  6.4M &   1.0M &  1.3k          &   14  &  3.1  &  160.8 &  2.4 &  163.3 \\ 
4096        & 50.9M &   8.4M & 11.5k          &  n/a  &  n/a  &   n/a  &  n/a &    n/a \\
\hline
\hline
\multicolumn{9}{|c|}{\textbf{3 levels}} \\
\hline
  64/8      &  0.8M &   0.1M &  0.1k/18       &   14  &  3.3  & 121.0  &  1.6 &  122.7 \\
 512/64     &  6.4M &   1.0M &  1.3k/295      &   17  &  4.2  & 166.9  &  2.4 &  169.3 \\
4096/512    & 50.9M &   8.4M & 11.5k/2930     &   18  &  4.4  & 221.7  &  5.5 &  227.3 \\
32768/128   & 405.0M & 69.1M & 95.2k/664      &   20  &  4.8  & 940.3  & 23.6 &  963.9 \\
\hline
\hline
\multicolumn{9}{|c|}{\textbf{4 levels}} \\
\hline
 512/64/8   &  6.4M &   1.0M &  1.3k/295/23   &   22  &  6.9  & 175.3  &  3.1 &  178.4 \\
4096/512/64 & 50.9M &   8.4M & 11.5k/2930/380 &   31  & 12.2  & 289.5  &  7.9 &  297.5 \\
32768/512/8 & 405.0M & 69.1M & 95.2k/2921/23  &   30  & 10.6  & 723.1  & 40.9 &  764.0 \\ 
\hline
\end{tabular}
\end{center}
\caption{\label{tab:cube_scaling_c1_adapt} 
Weak scaling for the cube problem with homogeneous material, \emph{adaptive multilevel BDDC}.
}
\end{table}

\begin{table}[tbph]
\begin{center}
\begin{tabular}
[c]{|c|ccc|cc|ccc|}
\hline
$N$ & 
\multirow{2}{*}{$n$} &
\multirow{2}{*}{$n_{\Gamma}$} &
$n_{f}$ &
\multirow{2}{*}{its.} & 
\multirow{2}{*}{cond.} & 
\multicolumn{3}{c|}{time (sec)} \\
$\ell=1(/2/3)$ &        
&
&        
$\ell=1(/2/3)$ & 
&    
&      
set-up & 
PCG & 
solve \\
\hline
\hline
\multicolumn{9}{|c|}{\textbf{2 levels}} \\
\hline
  8         &  0.1M &   9.5k &  12            &   582 & 236k &  4.0 &  59.4  &   63.4 \\
  64        &  0.8M &   0.1M &  0.1k          &  1611 & 233k &  4.7 & 171.9  &  176.6 \\
 512        &  6.4M &   1.0M &  1.3k          &  2195 & 240k &  9.5 & 340.4  &  350.0 \\
4096        & 50.9M &   8.4M & 11.5k          &   n/a & n/a  &  n/a &  n/a   &    n/a \\
\hline
\hline
\multicolumn{9}{|c|}{\textbf{3 levels}} \\
\hline
  64/8      &  0.8M &   0.1M &  0.1k/18       &  2218 & 239k &  4.7 &  234.1 &  238.8 \\
 512/64     &  6.4M &   1.0M &  1.3k/295      &  2830 & 250k &  5.5 &  328.2 &  333.7 \\
4096/512    & 50.9M &   8.4M & 11.5k/2930     &  4636 & 587k & 19.3 & 1096.2 & 1115.5 \\
32768/128   & 405.0M & 69.1M & 95.2k/664      &  6914 & 737k & 155.0 & 3820.8 & 3975.8 \\
\hline
\hline
\multicolumn{9}{|c|}{\textbf{4 levels}} \\
\hline
 512/64/8   &  6.4M &   1.0M &  1.3k/295/23   & 3771 &  729k & 5.4 &  434.4  & 439.8  \\
4096/512/64 & 50.9M &   8.4M & 11.5k/2930/380 & 8548 & 1860k & 9.3 & 1502.3  & 1511.6 \\
32768/512/8 & 405.0M & 69.1M & 95.2k/2921/23  & 9532 & 2362k & 160.2 & 5096.6 & 5256.8 \\ 
\hline
\end{tabular}
\end{center}
\caption{\label{tab:cube_scaling_c1e6_nonadapt} 
Weak scaling for the cube problem with jump in coefficients $E_2/E_1 = 10^6$, \emph{non-adaptive multilevel BDDC}.
}
\end{table}

\begin{table}[tbph]
\begin{center}
\begin{tabular}
[c]{|c|ccc|cc|ccc|}
\hline
$N$ & 
\multirow{2}{*}{$n$} &
\multirow{2}{*}{$n_{\Gamma}$} &
$n_{f}$ &
\multirow{2}{*}{its.} & 
\multirow{2}{*}{cond.} & 
\multicolumn{3}{c|}{time (sec)} \\
$\ell=1(/2/3)$ &        
&
&        
$\ell=1(/2/3)$ & 
&    
&      
set-up & 
PCG & 
solve \\
\hline
\hline
\multicolumn{9}{|c|}{\textbf{2 levels}} \\
\hline
  8         &  0.1M &  9.5k  &  12            &  119  & 1951 &   34.1 &  12.3   &   46.5  \\
  64        &  0.8M &  0.1M  &  0.1k          &   76  &  102 &   96.0 &   8.1   &  104.1  \\
 512        &  6.4M &  1.0M  &  1.3k          &   58  &   55 &  164.2 &   8.9   &  173.2  \\
4096        & 50.9M &  8.4M  & 11.5k          &   n/a &  n/a &  n/a   &   n/a   &    n/a  \\
\hline
\hline
\multicolumn{9}{|c|}{\textbf{3 levels}} \\
\hline
  64/8      &  0.8M &  0.1M  &  0.1k/18       & 457   &  48k &   96.7 &  48.0   &  144.7  \\ 
 512/64     &  6.4M &  1.0M  &  1.3k/295      & 82    & 0.1k &  165.7 &  10.2   &  175.9  \\ 
4096/512    & 50.9M &  8.4M  & 11.5k/2930     & 282   & 165k &  238.7 &  74.1   &  312.9  \\ 
32768/128   & 405.0M & 69.1M & 95.2k/664      & 270   &  24k  & 909.4 & 297.6   &  1207.0 \\
\hline
\hline
\multicolumn{9}{|c|}{\textbf{4 levels}} \\
\hline
 512/64/8   &  6.4M &   1.0M &  1.3k/295/23   & 554   &  63k &  169.5 &  68.3   &  273.7  \\
4096/512/64 & 50.9M &   8.4M & 11.5k/2930/380 & 3392  & 671k &  299.3 & 800.1   & 1099.4  \\
32768/512/8 & 405.0M & 69.1M & 95.2k/2921/23  & 3762  & 10495k & 697.6 & 4925.1 & 5622.7 \\ 
\hline
\end{tabular}
\end{center}
\caption{\label{tab:cube_scaling_c1e6_adapt} 
Weak scaling for the cube problem with jump in coefficients $E_2/E_1 = 10^6$, \emph{adaptive multilevel BDDC}.
}
\end{table}

\begin{table}[tbph]
\begin{center}
\begin{tabular}
[c]{|c|ccc|ccc|c|cc|}
\hline
\multirow{2}{*}{$\tau ^2$} & 
\multicolumn{3}{c|}{1st level (11\,520 pairs)} &
\multicolumn{3}{c|}{2nd level (144 pairs)} &
\multirow{2}{*}{$\widetilde{\omega}$} &
\multirow{2}{*}{its.} &
\multirow{2}{*}{cond.} \\
 &
ad. cstrs. &
cstrs./pair &
$\widetilde{\omega} _1$ &
ad. cstrs. &
cstrs./pair &
$\widetilde{\omega} _2$ &
& 
&
\\
\hline
\hline
25.0 &   120   & 0.01 & 4.45 &   84 &  0.58 & 4.37 & 19.53 & 21 & 7.18 \\
16.0 &  2220   & 0.19 & 2.70 &  132 &  0.92 & 3.33 &  8.99 & 18 & 4.77 \\
9.0 &  2220   & 0.19 & 2.70 &  228 &  1.58 & 2.92 &  7.89 & 18 & 4.72 \\
4.0 & 15\,660 & 1.36 & 1.99 & 1116 &  7.75 & 1.98 &  3.93 & 13 & 2.77 \\
2.25 & 69\,960 & 6.07 & 1.42 & 1440 & 10.00 & 2.49 &  3.55 & 14 & 3.25 \\
\hline
\end{tabular}
\end{center}
\caption{\label{tab:adaptive_constraints_on_cube} 
Comparison of adaptively selected constraints for different target condition number $\tau ^2$; 
`ad. cstrs.' -- number of added adaptive constraints, 
`cstrs./pair' -- average number of constraints added for one pair,
$\widetilde{\omega} = \widetilde{\omega} _1 \widetilde{\omega} _2$ -- the condition number indicator from 
(\ref{eq:cond-ind-mlevel});
\emph{adaptive 3-level BDDC}.
}
\end{table}

\subsection{Elasticity in a cube with variable size of regions of jumps in coefficients}

\begin{figure}[ptbh]
\begin{center}
\begin{tabular}{c}
\includegraphics[width=0.3\textwidth]{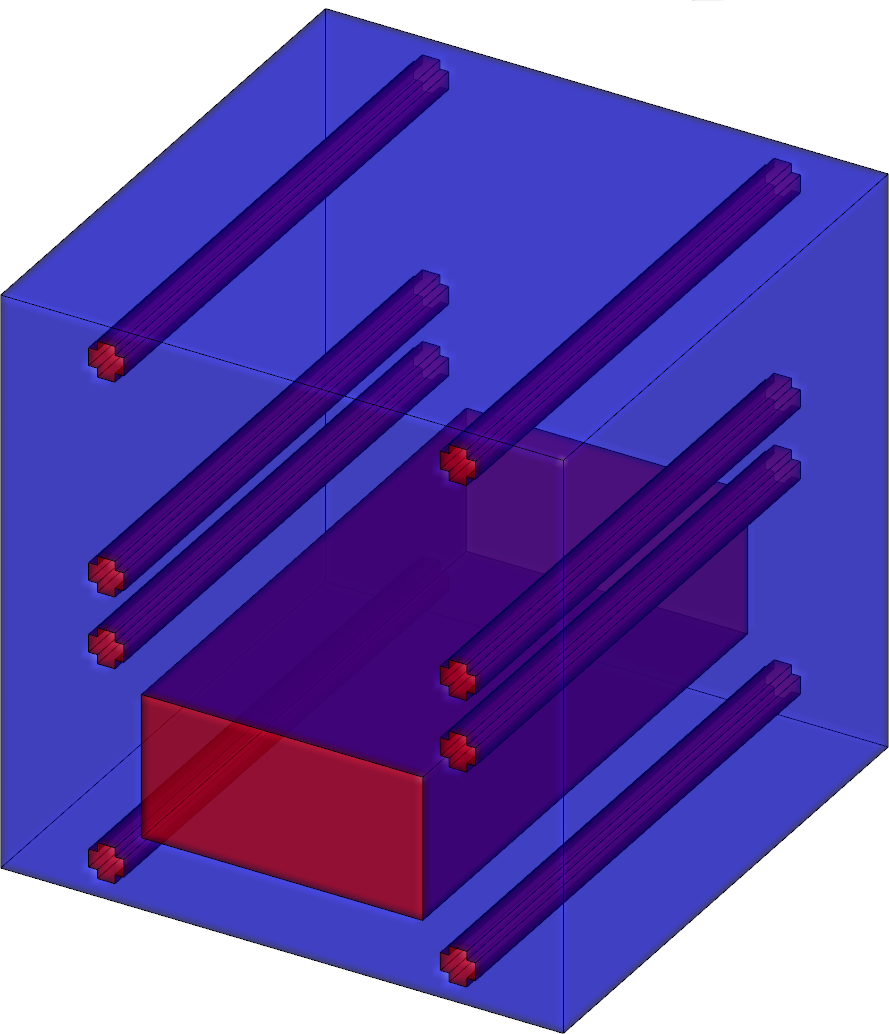}
\hskip 0.01\textwidth
\includegraphics[width=0.3\textwidth]{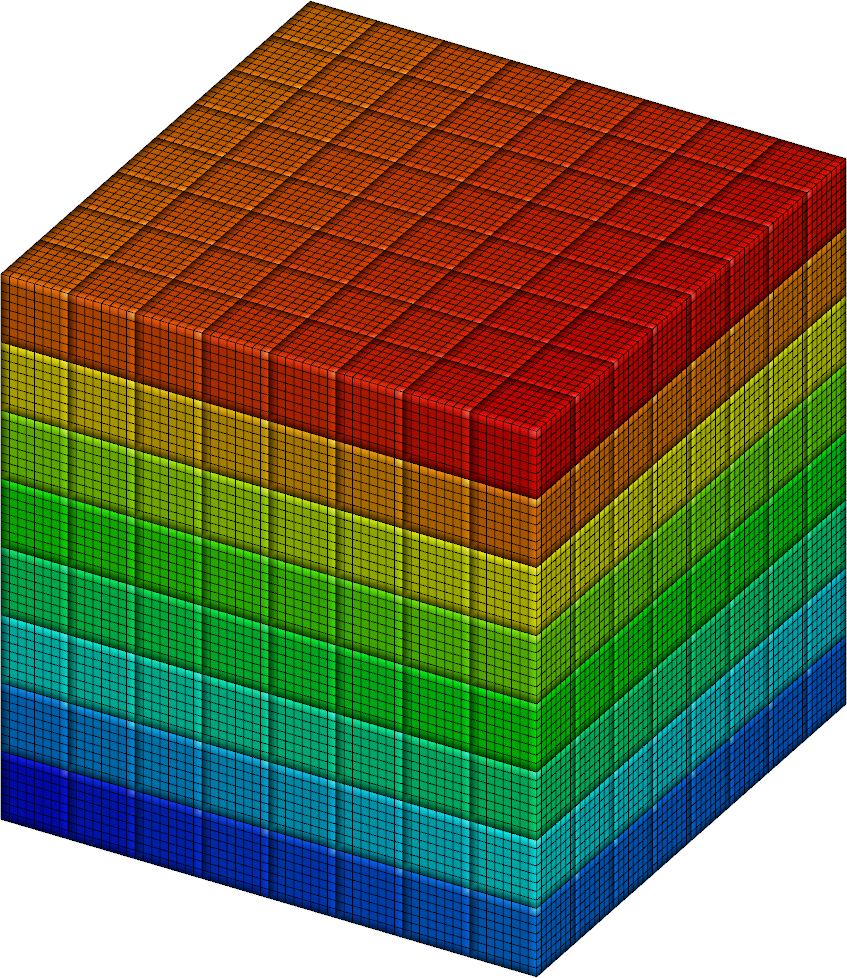}
\includegraphics[width=0.35\textwidth]{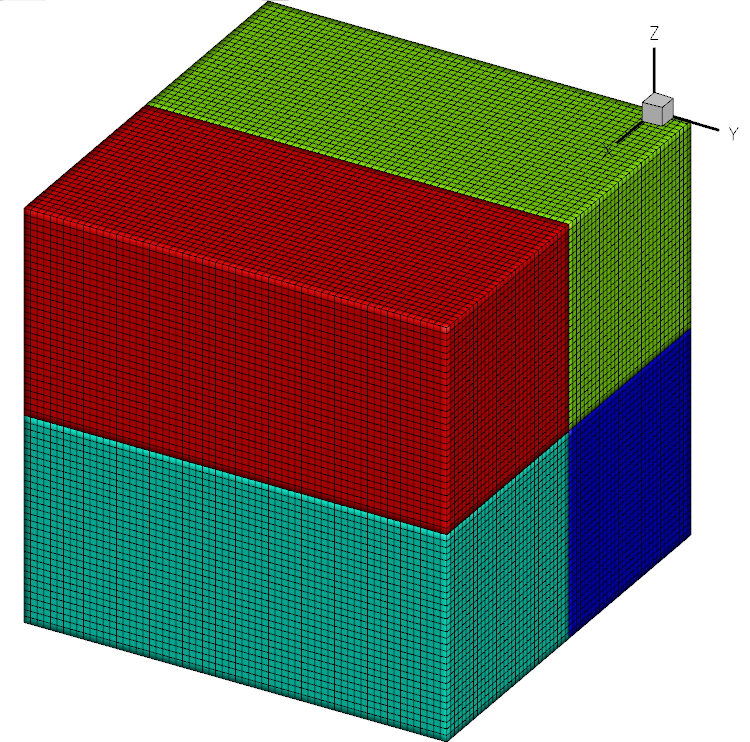}
\end{tabular}
\end{center}
\caption{\label{fig:multi2}
Cube with variable size of regions of jumps in coefficients: distribution of material with stiff bars ($E_2=2.1\cdot10^{11}$, $\nu_2= 0.3$), 
and soft outer material ($E_1=10^{6}$, $\nu_1=0.45$) (left),
mesh consisting of 823k degrees of freedom distributed into 512 substructures with 1344 faces on the first decomposition level (centre), 
and 4 substructures with 4 faces on the second decomposition level (right). 
Reproduced from~\cite{Sousedik-2010-AMB-thesis}.}
\end{figure}

The performance of the \emph{Adaptive-Multilevel BDDC} method in the presence of
jumps in material coefficients has been tested on a cube designed similarly
as the problem above, with material properties
$E_1=10^{6}$, $\nu_1=0.45$, and $E_2=2.1\cdot10^{11}$, $\nu_2= 0.3$. 
However, the stiff bars now vary in size, and while the thin bars create numerical difficulties on the first level, 
the large bar creates a jump in the decomposition on the second level, see~Fig.~\ref{fig:multi2}.
The computational mesh consists of 823k degrees of freedom and it is distributed into 512 substructures with 1344 faces on the
first decomposition level, and into 4 substructures with 4 faces on the second decomposition level (see~Fig.~\ref{fig:multi2}).

First, we present results by our serial implementation in Matlab, published initially in the thesis \cite{Sousedik-2010-AMB-thesis}.
We include them here along the parallel results to make this study of \emph{Adaptive-Multilevel BDDC} more self-contained. 
Comparing the results in Tabs.~\ref{tab:multi2-2lev-nonadaptive} and~\ref{tab:multi2-2lev-adaptive} we see that 
a~relatively small number of (additional) constraints leads to a considerable
decrease in number of iterations of the 2-level method. 
In these tables, 
$Nc$ denotes number of constraints, 
`c',`c+e', `c+e+f' denote combinations of constraints at corners, and arithmetic averages at edges and faces,
`3eigv' corresponds to using three adaptive constraints on faces instead of the three arithmetic averages,
$\tau$ denotes the target condition number from Algorithm~\ref{alg:setup_adaptive},
$\widetilde{\omega}$ is the indicator of the condition number from (\ref{eq:cond-ind-mlevel}),
`cond.' denotes estimated condition number, and 
`its.' the number of PCG iterations.

When the non-adaptive 2-level is replaced by the 3-level method (Tabs.~\ref{tab:multi2-2lev-nonadaptive} and \ref{tab:multi2-3lev-nonadaptive}), 
the condition number estimate as well as the number of iterations grow,
in agreement with the estimate (\ref{eq:def-omega-mlevel}).
However, with the adaptive 3-level approach (Tab.~\ref{tab:multi2-3lev-adaptive}) we were able to achieve nearly the same
convergence properties for small $\tau $ as in the adaptive 2-level method (Tab.~\ref{tab:multi2-2lev-adaptive}).

\begin{table}[ptbh]
\begin{center}
\begin{tabular}{|c|rrr|}
\hline
constraint       & $Nc$     & cond.      & its.      \\ \hline
c                & $2\,163$ & $312\,371$ & $>3\,000$ \\
c+e              & $5\,691$ & $45\,849$  & $1\,521$  \\
e+e+f            & $9\,723$ & $16\,384$  & $916$     \\
c+e+f$\,$(3eigv) & $9\,723$ & $3\,848$   & $367$     \\
\hline
\end{tabular}
\end{center}
\caption{\label{tab:multi2-2lev-nonadaptive}
Results for the cube with variable size of regions of jumps in coefficients (Fig.~\ref{fig:multi2}) obtained
using the \emph{non-adaptive 2-level BDDC} method. 
Reproduced from \cite{Sousedik-2010-AMB-thesis}.}
\end{table}

\begin{table}[ptbh]
\begin{center}
\begin{tabular}{|r|rrrr|}
\hline
$\tau$                    & $Nc$      & $\widetilde{\omega}$ & cond.        & its. \\
\hline
$\infty${\small {(=c+e)}} & $5\,691$  & ${\mathcal O}(10^{4})$          & $45\,848.60$ & $1\,521$ \\
$10\,000$                 & $5\,883$  & $8\,776.50$          & $5\,098.60$  & $441$ \\
$1\,000$                  & $6\,027$  & $5.33 $              & $9.92$       & $32$ \\
$10$                      & $6\,149$  & $6.25 $              & $6.66$       & $28$ \\
$5$                       & $9\,119$  & $<5 $                & $4.79$       & $24$ \\
$2$                       & $25\,009$ & $<2 $                & $2.92$       & $18$ \\
\hline
\end{tabular}
\end{center}
\caption{\label{tab:multi2-2lev-adaptive}
Results for the cube with variable size of regions of jumps in coefficients (Fig.~\ref{fig:multi2}) obtained
using the \emph{adaptive 2-level BDDC} method. 
Reproduced from \cite{Sousedik-2010-AMB-thesis}.}
\end{table}

\begin{table}[tbph]
\begin{center}
\begin{tabular}{|c|rrr|}
\hline
constraint       & $Nc$        & cond.       & its.       \\ \hline
c                & $2163/18$   & ${\mathcal O}(10^{7})$ & $>3\,000$ \\
c+e              & $5691/21$   & ${\mathcal O}(10^{6})$ & $>3\,000$  \\
c+e+f            & $9723/33$   & $461\,750$  & $1\,573$   \\
c+e+f$\,$(3eigv) & $9723/33$   & $125\,305$  & $981$      \\
\hline
\end{tabular}
\end{center}
\caption{\label{tab:multi2-3lev-nonadaptive} 
Results for the cube with variable size of regions of jumps in coefficients (Fig.~\ref{fig:multi2}) obtained
using the \emph{non-adaptive 3-level BDDC} method. 
Reproduced from \cite{Sousedik-2010-AMB-thesis}.}
\end{table}

\begin{table}[tbph]
\begin{center}
\begin{tabular}{|r|rrrr|}
\hline
$\tau ^2$                    & $Nc$           & $\widetilde{\omega}$ & cond.                  & its.       \\
\hline
$\infty${\small {(=c+e)}} & $5691+ 21$     & -                    & ${\mathcal O}(10^{6})$ & $ >3000$   \\
$10\,000$                 & $5883/28$      & $8776.50$            & $26\,874.40$           & $812$      \\
$1000$                    & $6027/34$      & $766.82$             & $1449.50$              & $145$      \\
$100$                     & $6027/53$      & $99.05$              & $100.89$               & $59$       \\
$10$                      & $6149/65$      & $7.93$               & $7.91$                 & $30$       \\
$5$                       & $9119/67$      & $<5 $                & $6.18$                 & $25$       \\
$2$                       & $25\,009/122$  & $<2 $                & $3.08$                 & $18$       \\
\hline
\end{tabular}
\end{center}
\caption{\label{tab:multi2-3lev-adaptive} 
Results for the cube with variable size of regions of jumps in coefficients (Fig.~\ref{fig:multi2}) obtained
using the \emph{adaptive 3-level BDDC} method. 
Reproduced from \cite{Sousedik-2010-AMB-thesis}.}
\end{table}

Next, we use this test problem to perform a~strong scaling test of our parallel implementation of adaptive multilevel BDDC method.
Since BDDCML supports assigning several subdomains to each processor, the division is kept constant with 512 subdomains on the basic level, 
and 4 subdomains on the second level (as in Fig.~\ref{fig:multi2}), and number of cores is varied.

\begin{figure}[tbph]
\begin{center}
\includegraphics[width=0.495\textwidth]{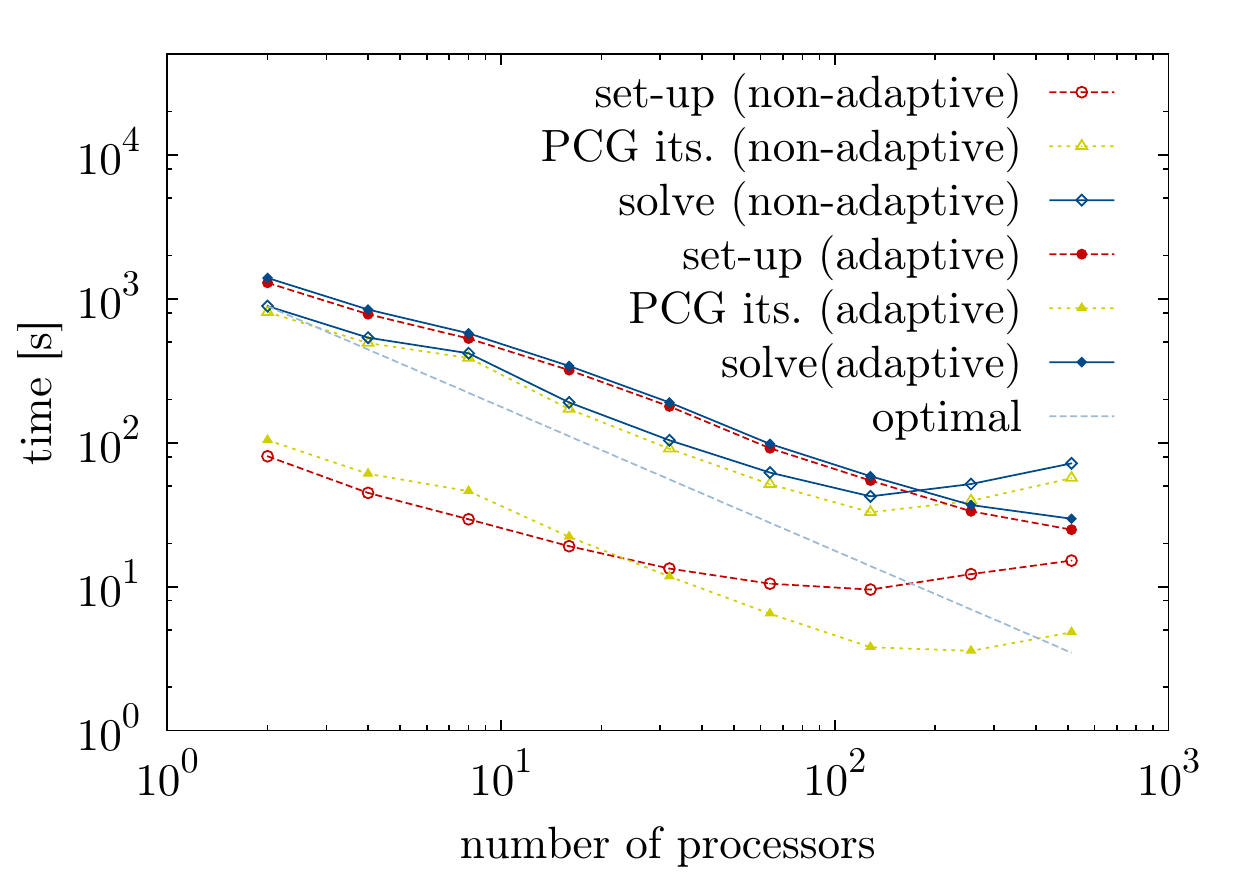}
\includegraphics[width=0.495\textwidth]{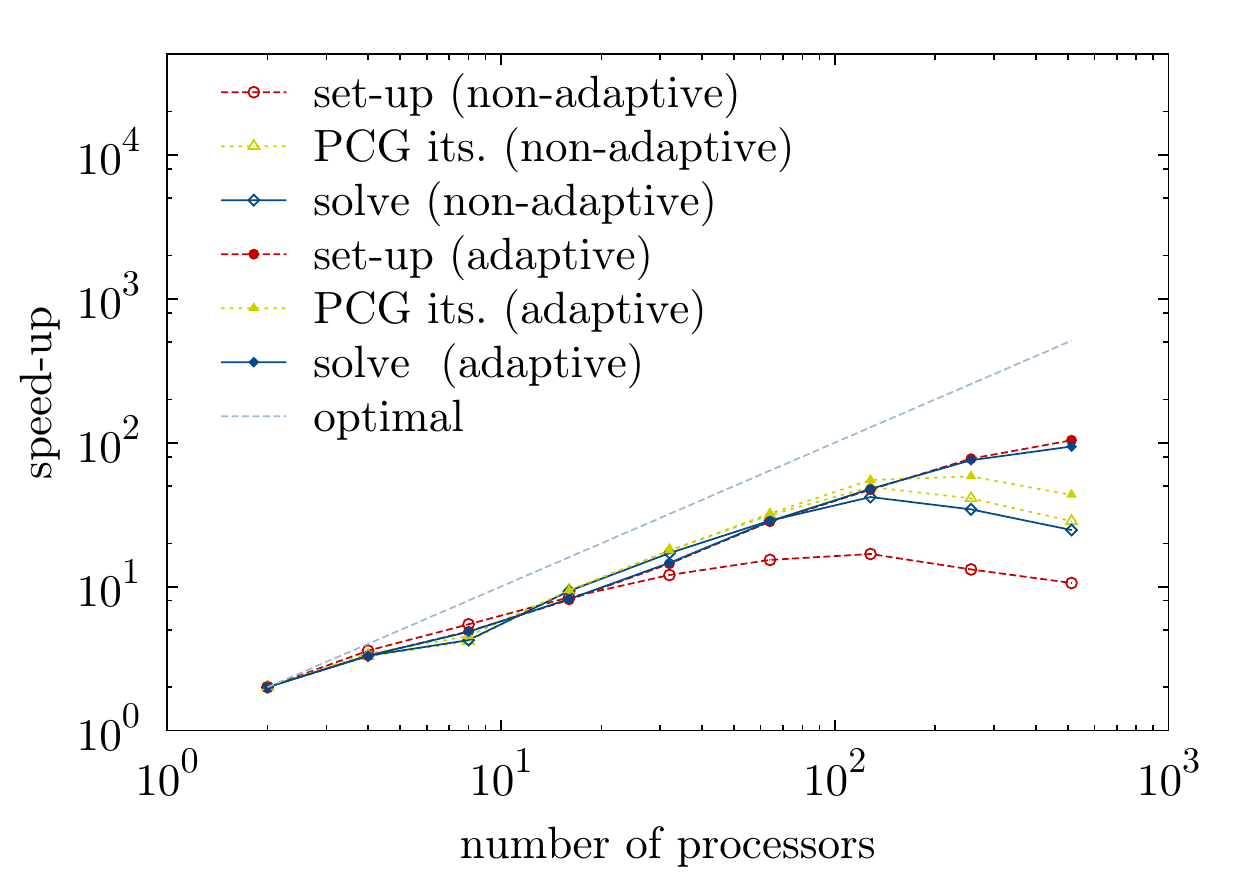}
\caption{\label{fig:timing_multi2_2level}
Strong scaling test for the cube with variable size of regions of jumps in coefficients (Fig.~\ref{fig:multi2}) containing 823k degrees of freedom,
on the first level divided into 512 subdomains with 1344 faces with arithmetic/adaptive constraints. 
Computational time (left) and speed-up (right) separately for set-up and PCG phases, and their sum (`solve'),
comparison of \emph{non-adaptive} (680 its.) and \emph{adaptive} (85 its.) \emph{parallel 2-level BDDC}. 
}
\end{center}
\end{figure}

\begin{figure}[tbph]
\begin{center}
\includegraphics[width=0.495\textwidth]{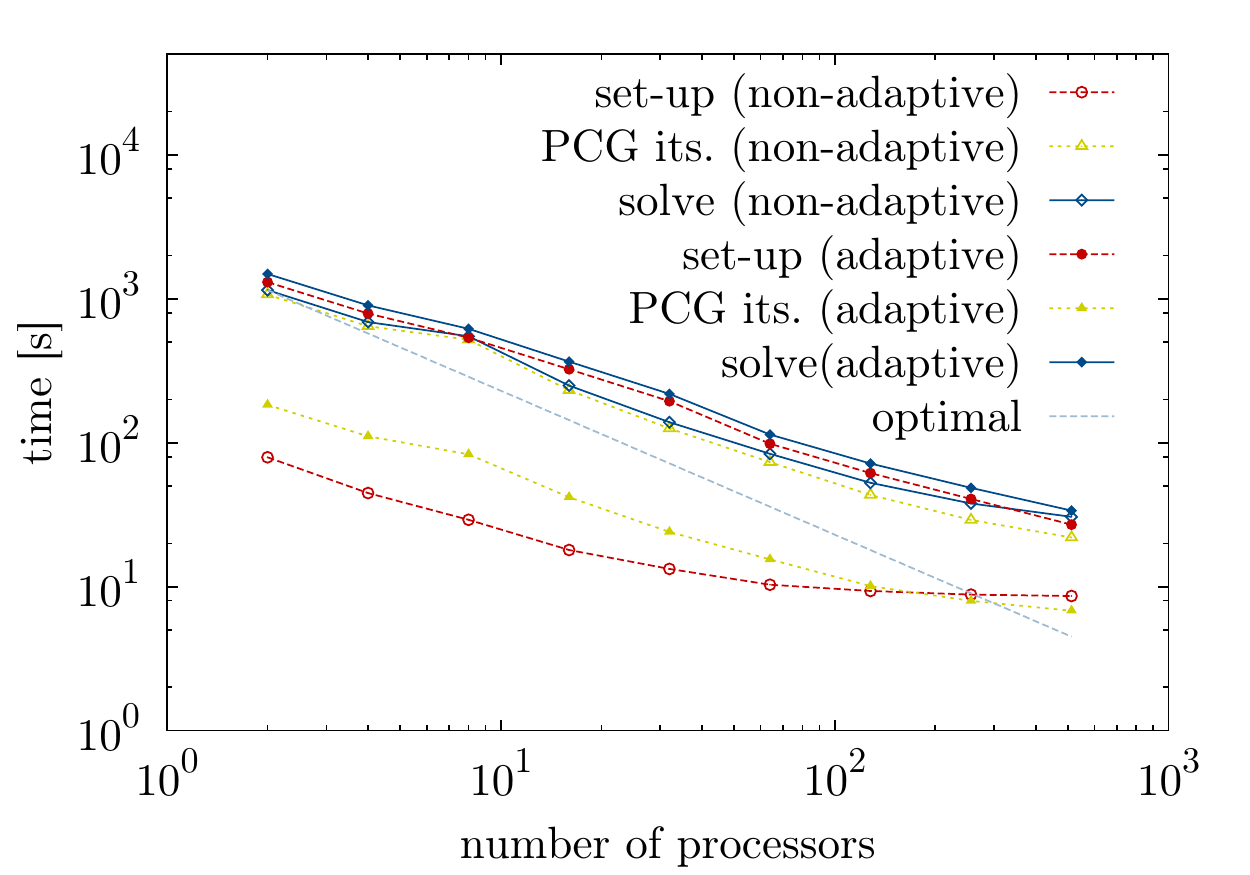}
\includegraphics[width=0.495\textwidth]{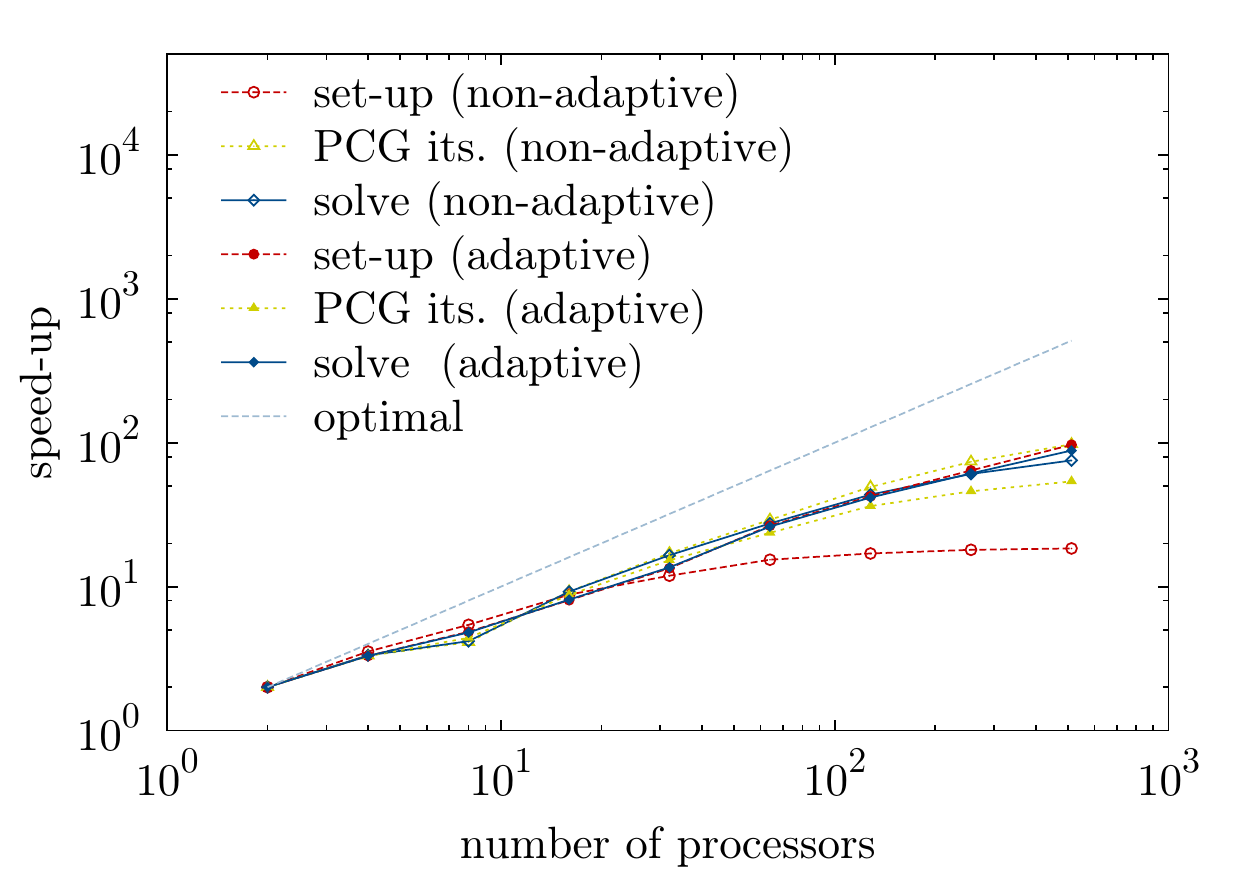}
\caption{\label{fig:timing_multi2_3level}
Strong scaling test for the cube with variable size of regions of jumps in coefficients (Fig.~\ref{fig:multi2}) containing 823k degrees of freedom,
on the first level divided into 512 subdomains with 1344 faces with arithmetic/adaptive constraints, and on the second level into 4 subdomains with 4 faces.
Computational time (left) and speed-up (right) separately for set-up and PCG phases, and their sum (`solve'),
comparison of \emph{non-adaptive} (894 its.) and \emph{adaptive} (150 its.) \emph{parallel 3-level BDDC}. 
}
\end{center}
\end{figure}

Figure~\ref{fig:timing_multi2_2level} presents parallel computational time and speed-up when this problem is solved by the parallel 
\emph{2-level BDDC} method, comparing efficiency of the non-adaptive and adaptive solver.
We report times and speed-ups independently for the set-up phase (including solution of eigenproblems for adaptive method),
the phase of PCG iterations, and their sum (`solve').
Figure~\ref{fig:timing_multi2_3level} then presents parallel computational time and speed-up for \emph{3-level BDDC} method.

We can see, that both phases of the solution are reasonably scalable. 
For large core counts, scalability worsens, as each core has only little work with subdomain problems and (the less scalable) solution  
of the coarse problem dominates the computation.
It is worth noting, that the overall fastest solution was delivered by the \emph{adaptive 2-level BDDC} method on 512 cores, 
while both other extensions of BDDC -- \emph{non-adaptive 3-level BDDC} and \emph{adaptive 3-level BDDC} -- were also considerably faster than the 
standard (\emph{non-adaptive 2-level}) BDDC method on large number of cores.

\subsection{Linear elasticity analysis of a mining reel}

\begin{figure}[ptbh]
\begin{center}
\includegraphics[width=0.47\textwidth]{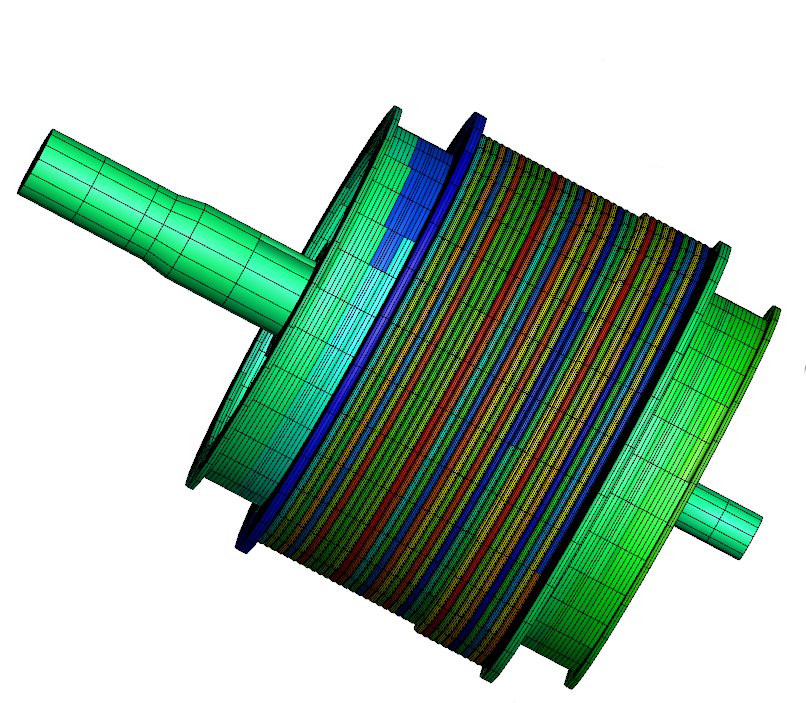} 
\includegraphics[width=0.47\textwidth]{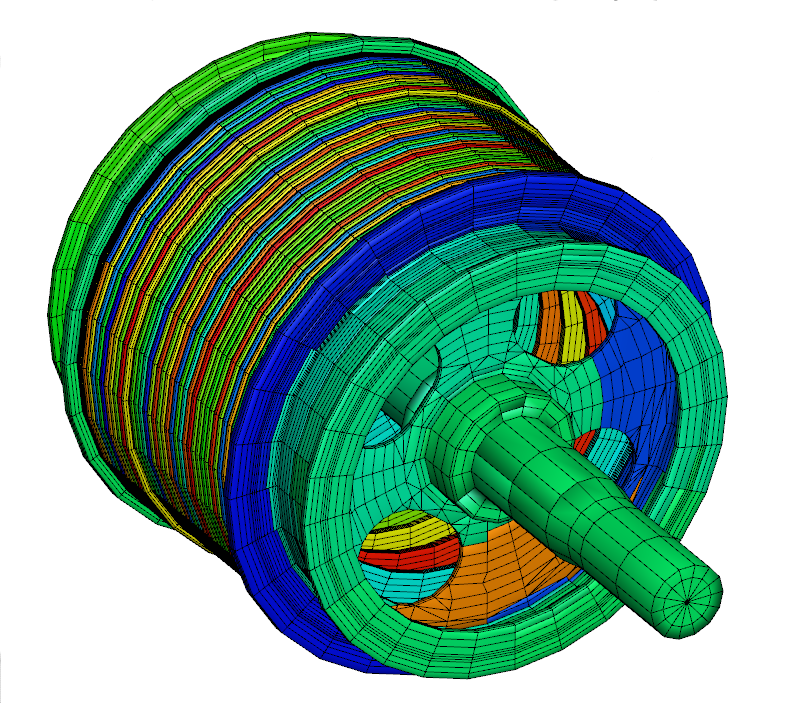} 
\end{center}
\caption{\label{fig:buben_1024}
Finite element discretization and substructuring of the mining reel
problem, consisting of 1.7M degrees of freedom, 
divided into 1024 subdomains with 3893 faces.
Data by courtesy of Jan Le\v{s}tina and Jaroslav Novotn\'{y}.
Reproduced from~\cite{Sousedik-2010-AMB-thesis}.
}
\end{figure}

The performance of the \emph{Adaptive-Multilevel BDDC} has been tested on 
an engineering problem of linear elasticity analysis of a mining reel. 
The problem was provided for testing by Jan Le\v{s}tina and Jaroslav Novotn\'{y}.
The computational mesh consists of 141k quadratic finite elements, 
579k nodes, and  approximately 1.7M degrees of freedom.
It was divided into 1024 subdomains with 3893 faces (see Fig.~\ref{fig:buben_1024}).

\begin{figure}[ptbh]
\begin{center}
\includegraphics[width=0.55\textwidth]{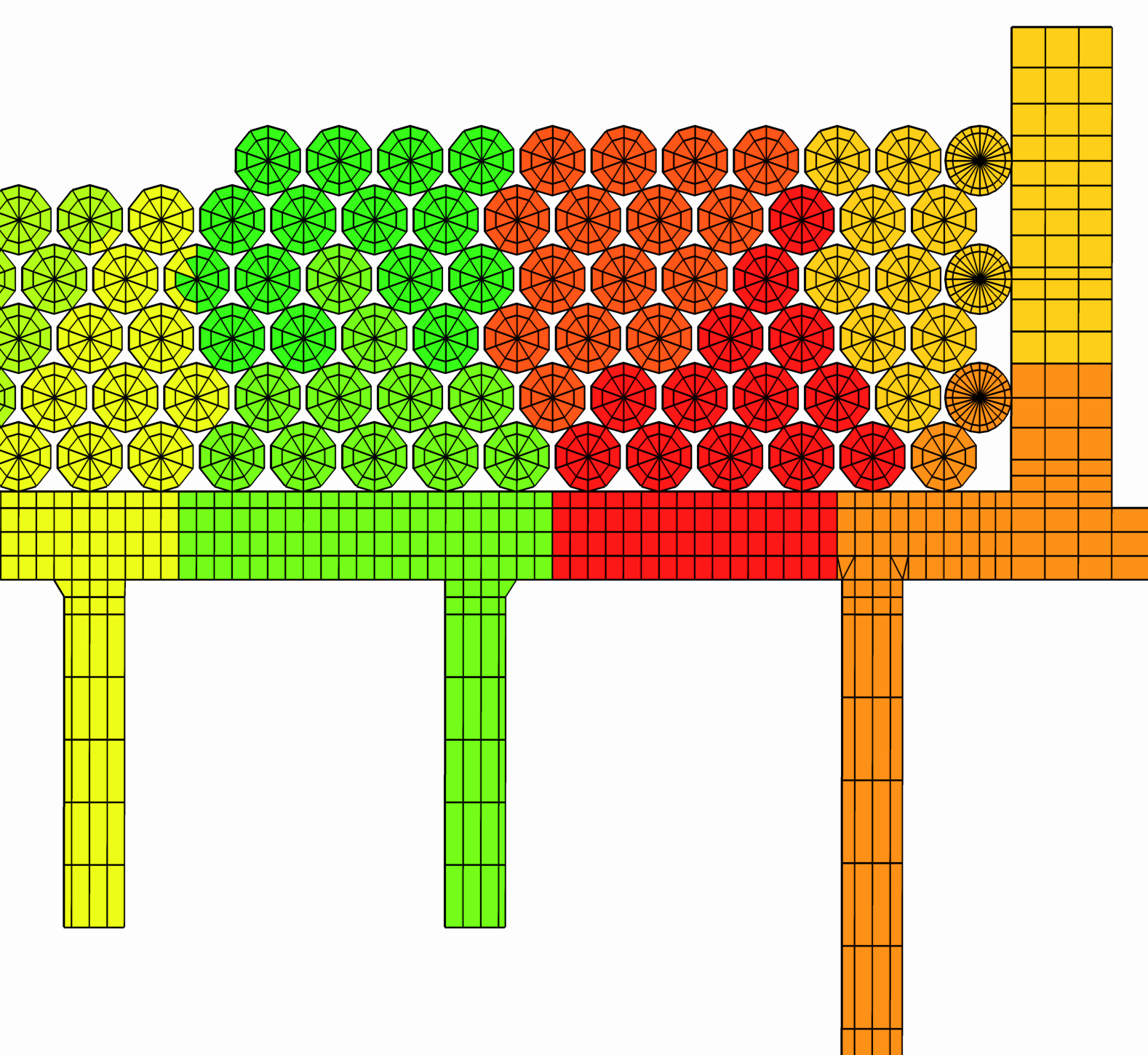} 
\hskip 5mm
\includegraphics[width=0.37\textwidth]{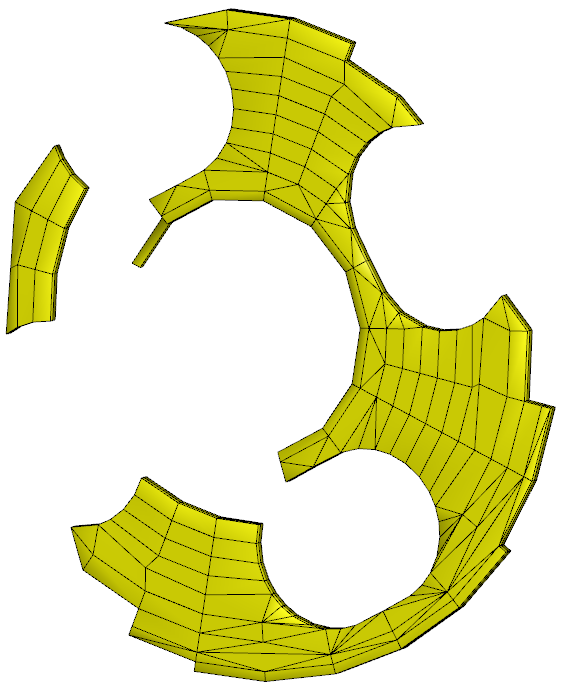} \\
\includegraphics[width=0.2\textwidth]{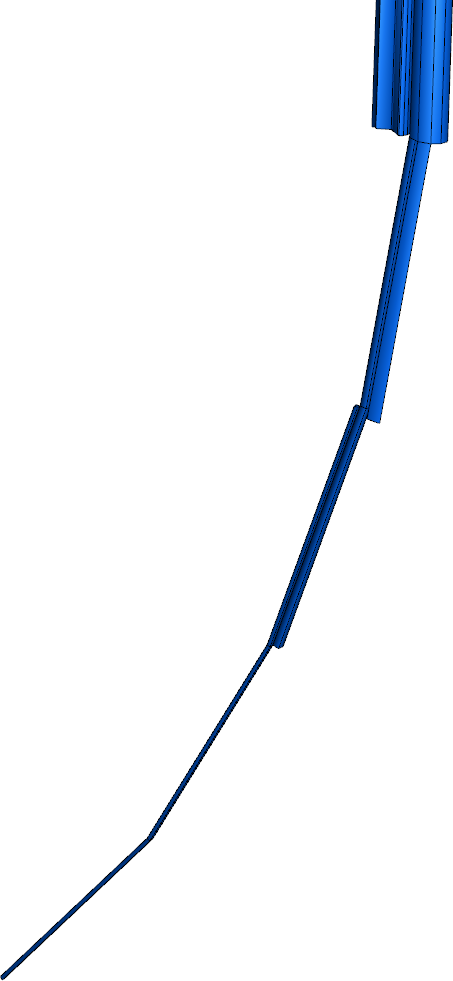} 
\includegraphics[width=0.65\textwidth]{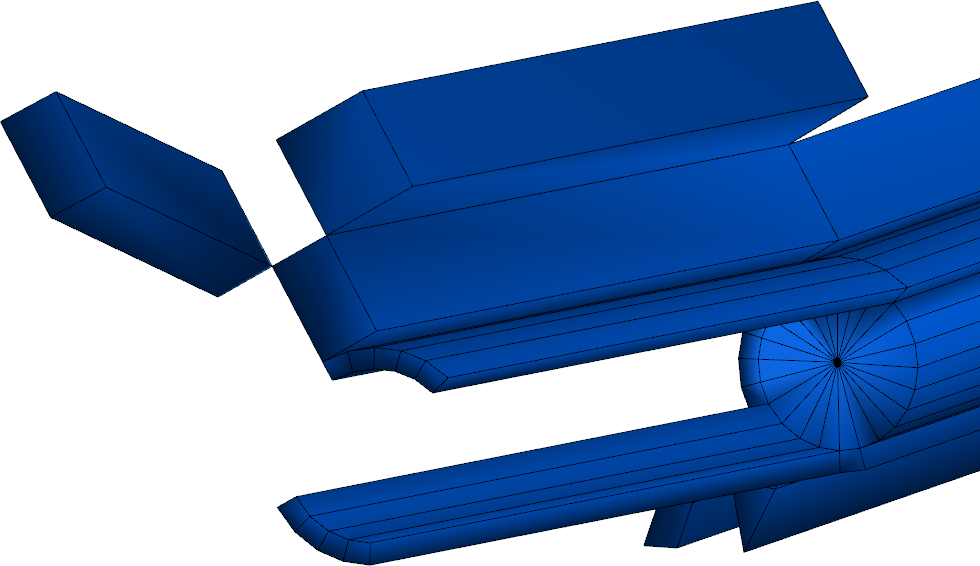} 
\end{center}
\caption{\label{fig:buben_1024_issues}
Examples of difficulties with computational mesh and its partitioning for the mining reel problem:
(i) rope modelled as axisymmetric rings connected at edges of three-dimensional elements (top left),
(ii) disconnected subdomains (top right),
(iii) elongated substructures (bottom left), 
(iv) spurious mechanisms within subdomains, such as elements connected to rest of the subdomain only at single node (bottom right).
Reproduced from~\cite{Sousedik-2010-AMB-thesis}.
}
\end{figure}

The problem presents a very challenging application for iterative solvers due to its very complicated geometry. 
It contains a steel rope, which is not modelled as a contact problem but just by a complicated mesh with elements connected through edges of 
three-dimensional elements (Fig.~\ref{fig:buben_1024_issues}).
Its automatic partitioning by METIS creates further problems such as thin elongated subdomains, 
disconnected subdomains, 
or subdomains with insufficiently coupled elements leading to `spurious mechanisms' inside subdomains.
See Fig.~\ref{fig:buben_1024_issues} for examples.

We first perform a~series of computations by our serial implementation in Matlab to study the effect of prescribed target condition number $\tau$ on convergence.
Comparing results by \emph{non-adaptive 2-level BDDC} (Tab.~\ref{tab:buben_1024-2lev-nonadaptive}) with \emph{adaptive 2-level BDDC} (Tab.~\ref{tab:buben_1024-2lev-adaptive}),
we see that the adaptive approach allows for a significant improvement in
the number of iterations. 

We can also see, that convergence of the adaptive two- and three-level method (Tables~\ref{tab:buben_1024-2lev-adaptive} and~\ref{tab:buben_1024-3lev-adaptive}) 
is nearly identical.
For the three-level method, automatic division into 32 subdomains was used on the second level.

We note that the observed approximate condition number computed
from the Lanczos sequence in PCG (`cond.') is larger than the target condition number $\tau$ for this problem. 
In \cite{Mandel-2007-ASF,Sousedik-2010-AMB-thesis}, it was shown that these two numbers match remarkably well for simpler problems, especially in 2D.
Despite of this difference, the algorithm still performs very well.

\begin{table}[tbph]
\begin{center}
\begin{tabular}{|c|rrr|}
\hline
constraint       & $Nc$      & cond.             & its.       \\ 
\hline
c+e              & $27\,183$ & -                 & $>2\,000$  \\
c+e+f            & $38\,868$ & $1.18\cdot10^{6}$ & $1\,303$   \\
c+e+f\ (3eigv)   & $38\,868$ & $72\,704.80 $     & $674$      \\
\hline
\end{tabular}
\end{center}
\caption{\label{tab:buben_1024-2lev-nonadaptive} Convergence of the \emph{non-adaptive 2-level BDDC method} with different constraints, mining reel problem.
Reproduced from \cite{Sousedik-2010-AMB-thesis}.
} 
\end{table}

\begin{table}[tbph]
\begin{center}
\begin{tabular}{|r|rrrr|}
\hline
$\tau$                    & $Nc$      & $\widetilde{\omega}$ & cond.    & its. \\
\hline
$\infty${\small {(=c+e)}} & $27\,183$ & $1.76\cdot10^{6}$    & -           & $ >2\,000 $  \\
$10000$                   & $28\,023$ & $9\,992.61 $         & $9\,538.18$ & $910$        \\
$5000$                    & $28\,727$ & $4\,934.62 $         & $4\,849.75$ & $673$        \\
$1000$                    & $32\,460$ & $999.90 $            & $2\,179.79$ & $391$        \\
$500$                     & $35\,017$ & $499.64 $            & $1\,277.59$ & $318$        \\
$100$                     & $42\,849$ & $99.89 $             & $840.74$    & $213$        \\
$50$                      & $46\,093$ & $49.98 $             & $784.49$    & $194$        \\
$10$                      & $59\,496$ & $<10$                & $321.20$    & $129$        \\
$5$                       & $69\,249$ & $<5$                 & $198.68$    & $91$         \\
$2$                       & $92\,467$ & $<2$                 & $91.24$     & $72$         \\
\hline
\end{tabular}
\end{center}
\caption{\label{tab:buben_1024-2lev-adaptive} Convergence of the \emph{adaptive 2-level BDDC method} with variable target condition number parameter $\tau$, 
mining reel problem. Lower $\tau$ corresponds to more constraints and better convergence.
Reproduced from \cite{Sousedik-2010-AMB-thesis}.
}
\end{table}

\begin{table}[tbph]
\begin{center}
\begin{tabular}{|r|rrrr|}
\hline
$\tau ^2$ & $Nc$             &  $\widetilde{\omega}$ & cond.         & its.  \\ \hline
$100$  & $42\,849+2\,378$ &  $99.89$          & $ 3\,567.02 $ & $382$ \\
$10$   & $59\,496+6\,419$ &  $<10$            & $320.82$      & $139$ \\
$5$    & $69\,249+8\,681$ &  $<5$             & $198.55$      & $98$  \\
\hline
\end{tabular}
\end{center}
\caption{\label{tab:buben_1024-3lev-adaptive} Convergence of the \emph{adaptive 3-level BDDC method} with variable target condition number parameter $\tau$, 
mining reel problem.
Reproduced from \cite{Sousedik-2010-AMB-thesis}.
}
\end{table}

\begin{figure}[tbph]
\begin{center}
\includegraphics[width=0.495\textwidth]{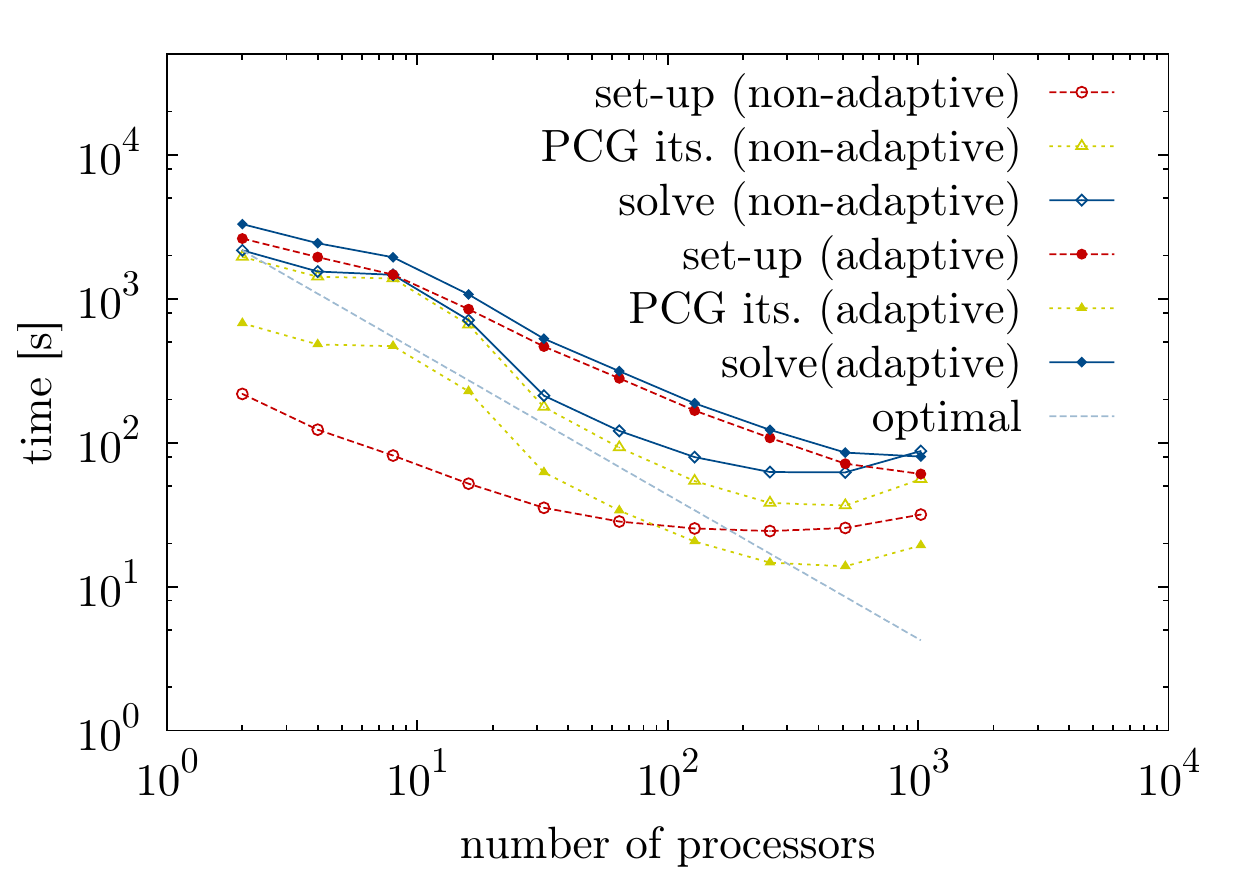}
\includegraphics[width=0.495\textwidth]{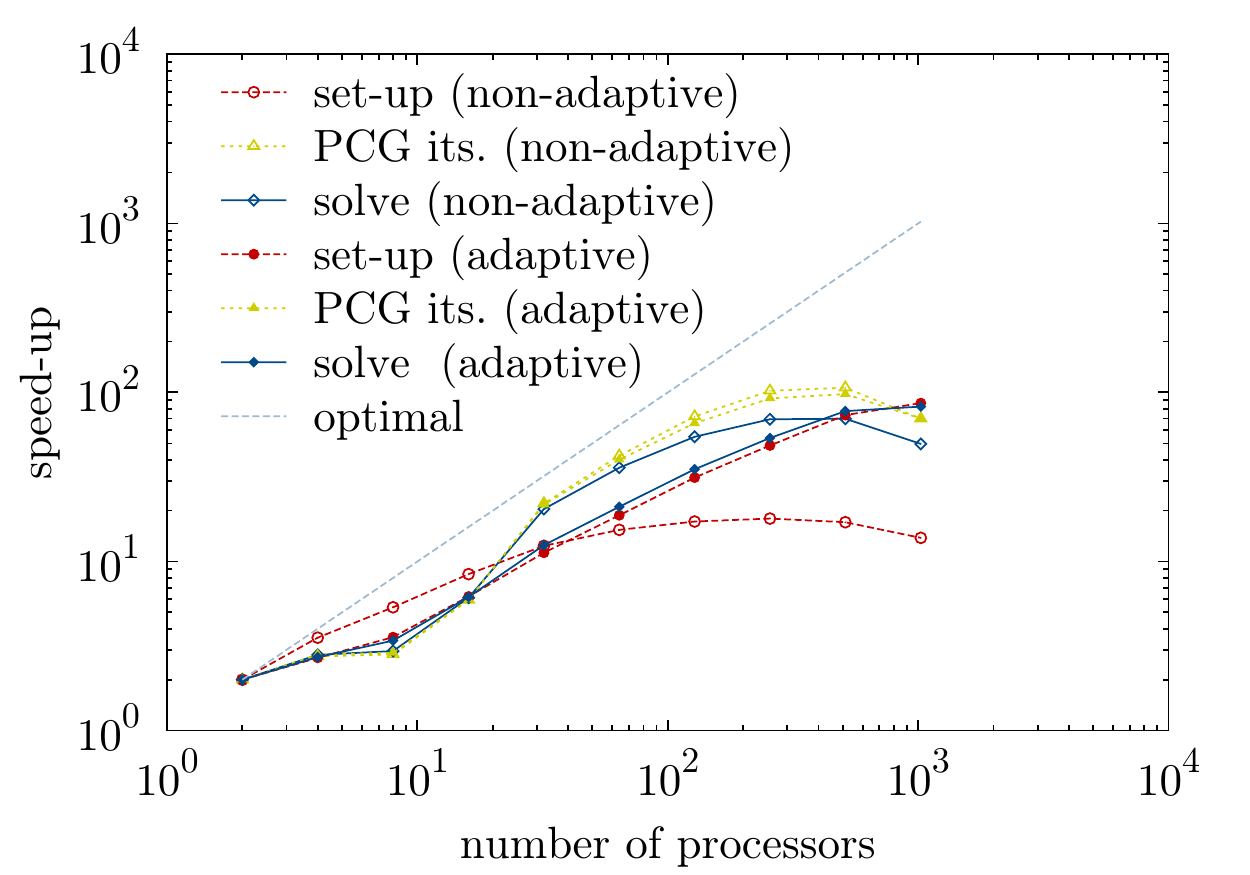}
\caption{\label{fig:timing_buben}
Strong scaling test for the mining reel problem containing 1.7M degrees of freedom and on the first level divided into 1024 subdomains with 3893 faces with arithmetic/adaptive constraints, computational time (left) and speed-up (right) separately for set-up and PCG phases, and their sum (`solve'),
comparison of \emph{non-adaptive} (610 its.) and \emph{adaptive} (200 its.) \emph{parallel 2-level BDDC}. 
}
\end{center}
\end{figure}

Next, we solved this problem with the parallel implementation of the algorithm.
In Fig.~\ref{fig:timing_buben}, we present results of a strong scaling test.

The \emph{non-adaptive 2-level BDDC} method requires 610 PCG iterations, while the \emph{adaptive 2-level BDDC} needs only 200 PCG iterations.
Nevertheless, this difference is only able to compensate the cost of solving the eigenproblems, 
and the adaptive method is advantageous with respect to computing time only for 1024 cores.

We can see, that the scaling is nearly optimal, 
with the deviation caused probably again by the small size of the problem compared to the core counts used in this experiment.
In the parallel case, the 3-level approach did not work well neither with nor without adaptive selection of constraints, 
requiring more than 5000 PCG iterations in both cases.
The slow convergence of the adaptive 3-level method is probably caused by the limit of ten adaptive constraints per face, 
which seems to be insufficient for this difficult problem -- in Matlab experiments, as many as 8681 adaptive constraints were generated 
among 32 subdomains on the second level in order to satisfy $\widetilde{\omega} \leq \tau ^2 = 5$ (Tab.~\ref{tab:buben_1024-3lev-adaptive}).

\subsection{Linear elasticity analysis of a geocomposite sample}

Finally, the algorithm is applied to a problem of elasticity analysis of a cubic geocomposite sample.
The sample was analyzed in \cite{Blaheta-2009-SDD} and provided by the authors for testing of our implementation.
The length of the edge of the cube is 75 mm, and 
the cube is composed of five distinct materials identified by means of computer tomography. 

\begin{figure}[ptbh]
\begin{center}
\includegraphics[width=0.47\textwidth]{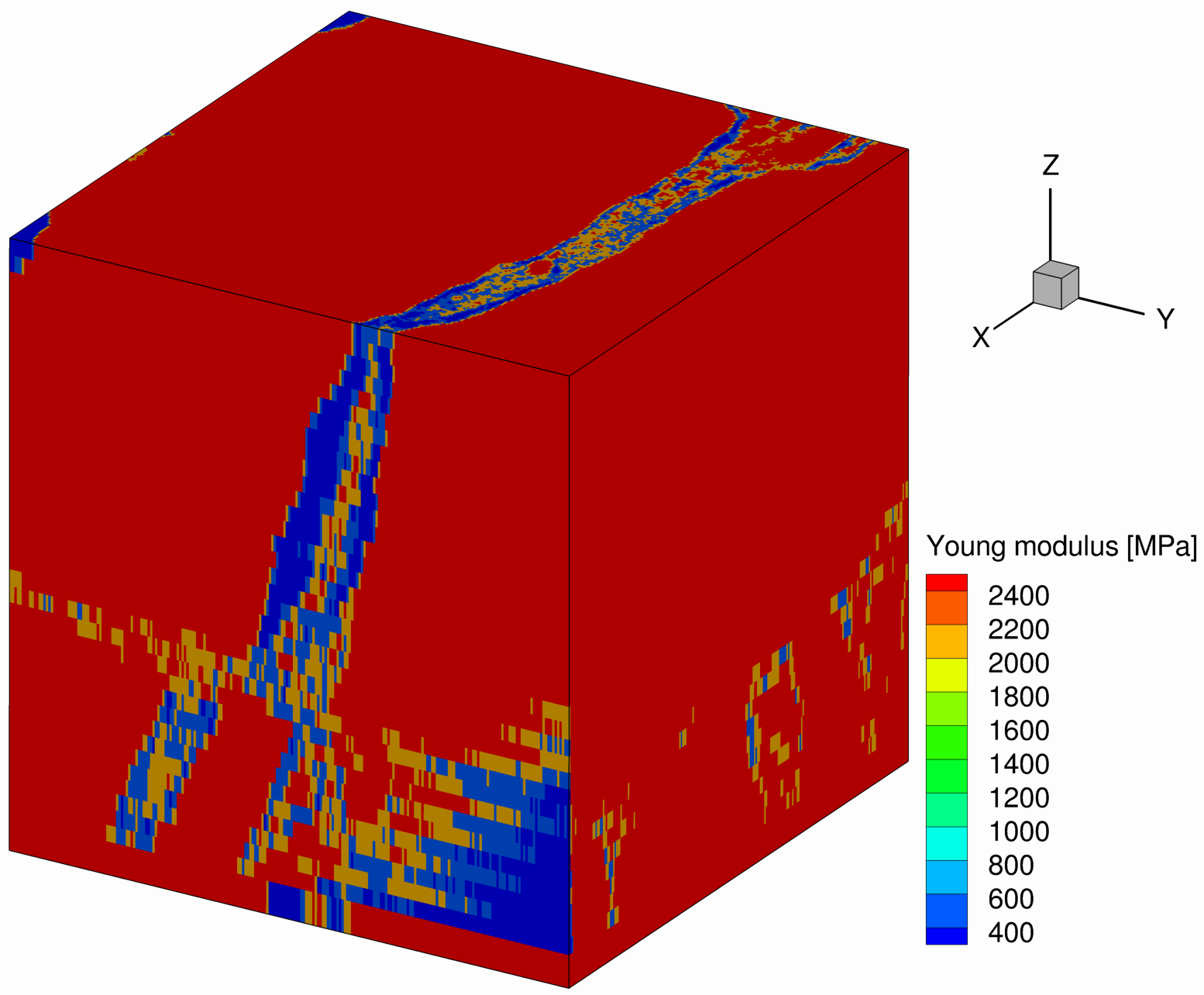} 
\includegraphics[width=0.47\textwidth]{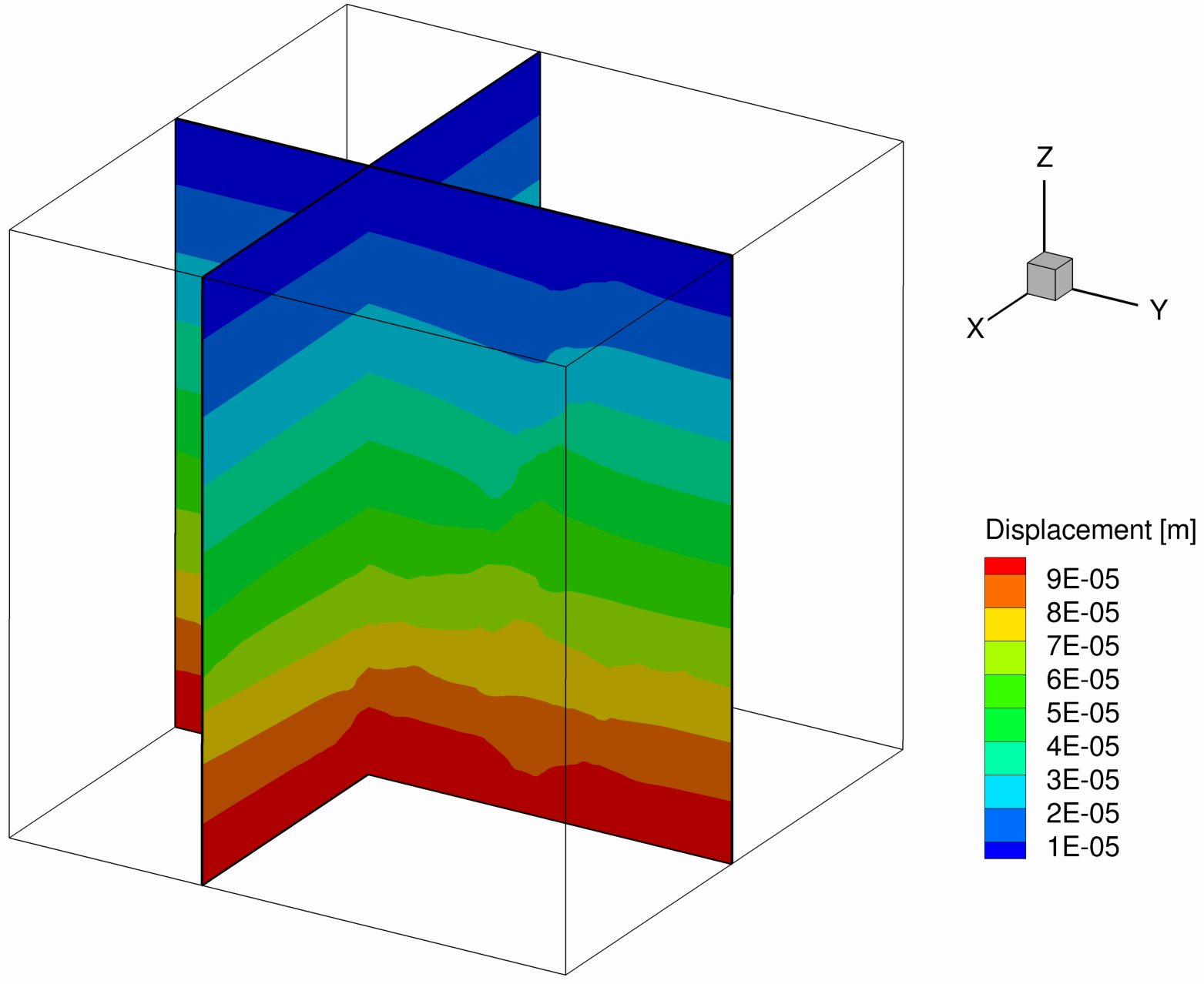} 
\end{center}
\caption{\label{fig:geoc}
Elasticity in a geocomposite sample; Young's modulus due to different materials (left), 
magnitude of displacement on slices (right).
Mesh contains 12 million linear tetrahedral elements and approx. 6 million degrees of freedom.
Data by courtesy of Radim Blaheta and Ji\v{r}\'{\i} Star\'{y}.
}
\end{figure}

Different material properties cause anisotropic response of the cube even for simple axial stretching in z direction 
(Fig. \ref{fig:geoc} right).
The problem is discretized using unstructured grid of about 12 million linear tetrahedral elements,
resulting in approximately 6 million degrees of freedom.
The mesh was divided into 1024 subdomains on the first and into 32 subdomains on the second level,
resulting in 5635 and 100 faces, respectively.

Table~\ref{tab:geoc_convergence} summarizes the number of iterations and estimated condition number of the preconditioned
operator for the combinations of 2- and 3-level BDDC with non-adaptive and adaptive selection of constraints.
For this problem, number of iterations was reduced to approximately one half by using adaptivity.
We can also see, that number of iterations grows considerably when going from 2 to 3 levels for this problem.

\begin{table}[tbph]
\begin{center}
\begin{tabular}
[c]{|c|c|ccc|cc|cc|}
\hline
\multirow{2}{*}{levels} &
$N$ & 
\multirow{2}{*}{$n$} &
\multirow{2}{*}{$n_{\Gamma}$} &
$n_{f}$ &
\multicolumn{2}{|c|}{non-adaptive} & 
\multicolumn{2}{|c|}{adaptive} \\
&
$\ell=1(/2)$ &        
&
&        
$\ell=1(/2)$ & 
its. & 
cond.  &
its. & 
cond.  \\
\hline
  2 &  1024    & 6.1M &  1.3M  &  5635        &    65 &    94.8  &   36  &    37.7 \\ 
  3 &  1024/32 & 6.1M &  1.3M  &  5635/100    &   214 &  2724.9  &  113  &  1792.2 \\
\hline
\end{tabular}
\end{center}
\caption{\label{tab:geoc_convergence} 
Number of iterations and condition number estimate for the geocomposite problem, 6 million degrees of freedom.
Comparison for two and three levels for non-adaptive and adaptive BDDC.
}
\end{table}

\begin{figure}[tbph]
\begin{center}
\includegraphics[width=0.495\textwidth]{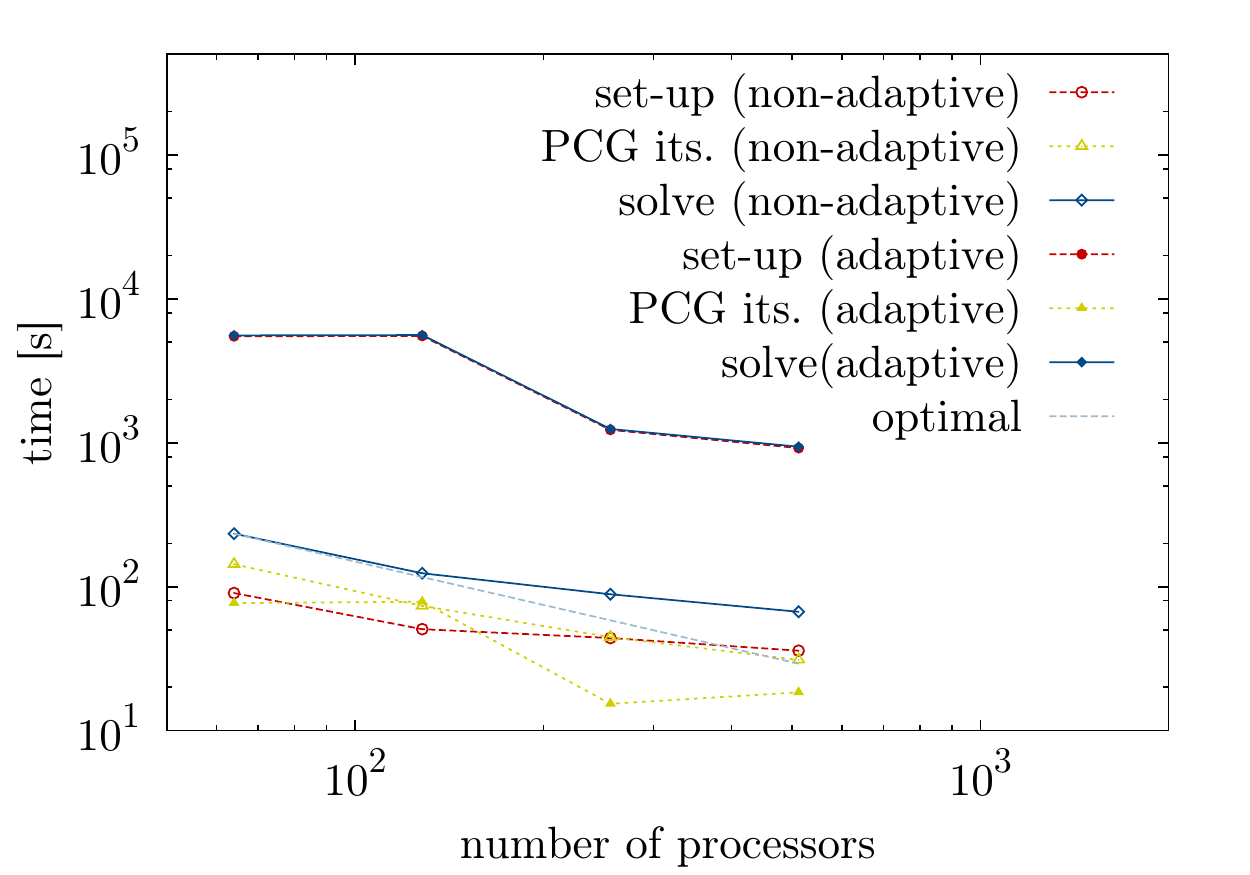}
\includegraphics[width=0.495\textwidth]{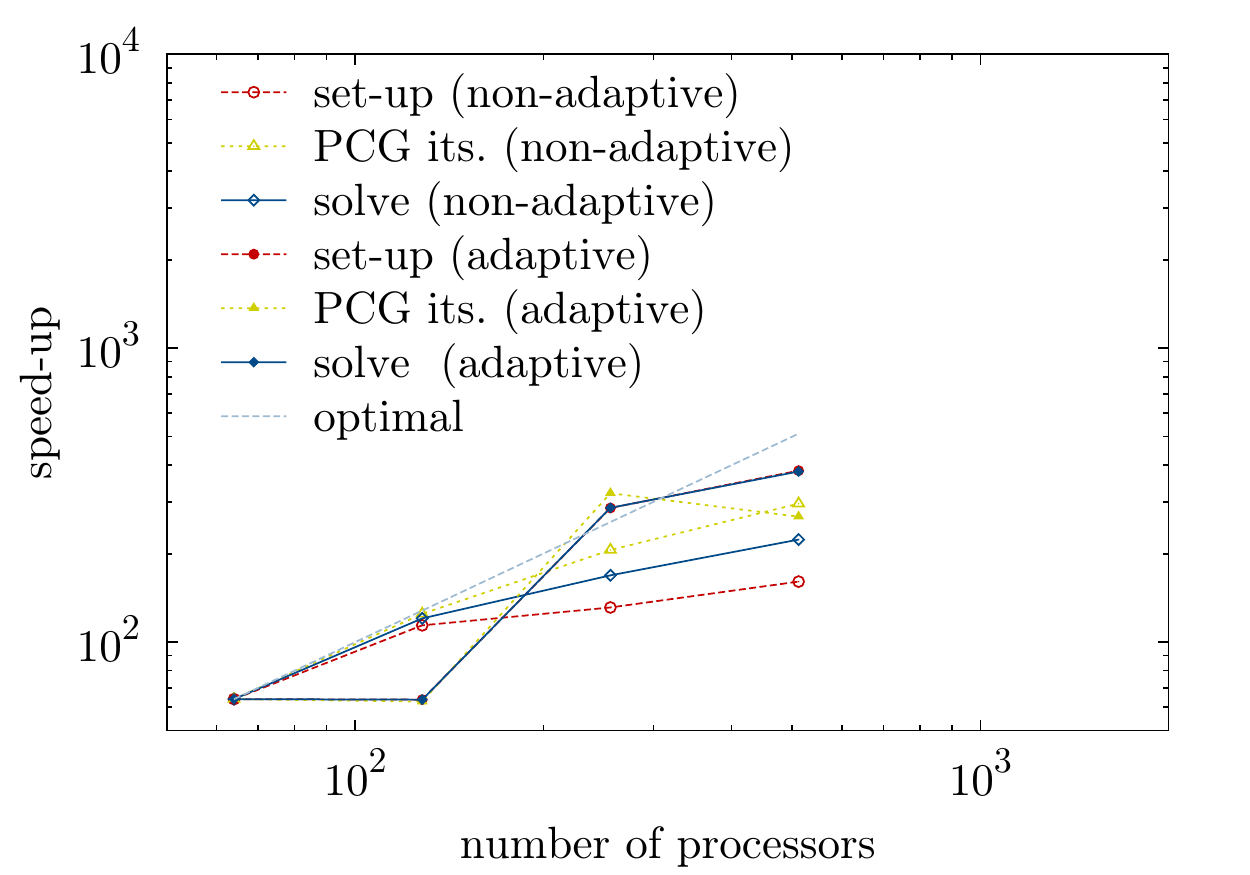}
\caption{\label{fig:timing_geoc1024_2level}
Strong scaling test for the geocomposite problem (Fig.~\ref{fig:geoc}) containing approx. 6 million degrees of freedom,
on the first level divided into 1024 subdomains with 5635 faces with arithmetic/adaptive constraints. 
Computational time (left) and speed-up (right) separately for set-up and PCG phases, and their sum (`solve'),
comparison of \emph{non-adaptive} (65 its.) and \emph{adaptive} (36 its.) \emph{parallel 2-level BDDC}. 
}
\end{center}
\end{figure}

\begin{figure}[tbph]
\begin{center}
\includegraphics[width=0.495\textwidth]{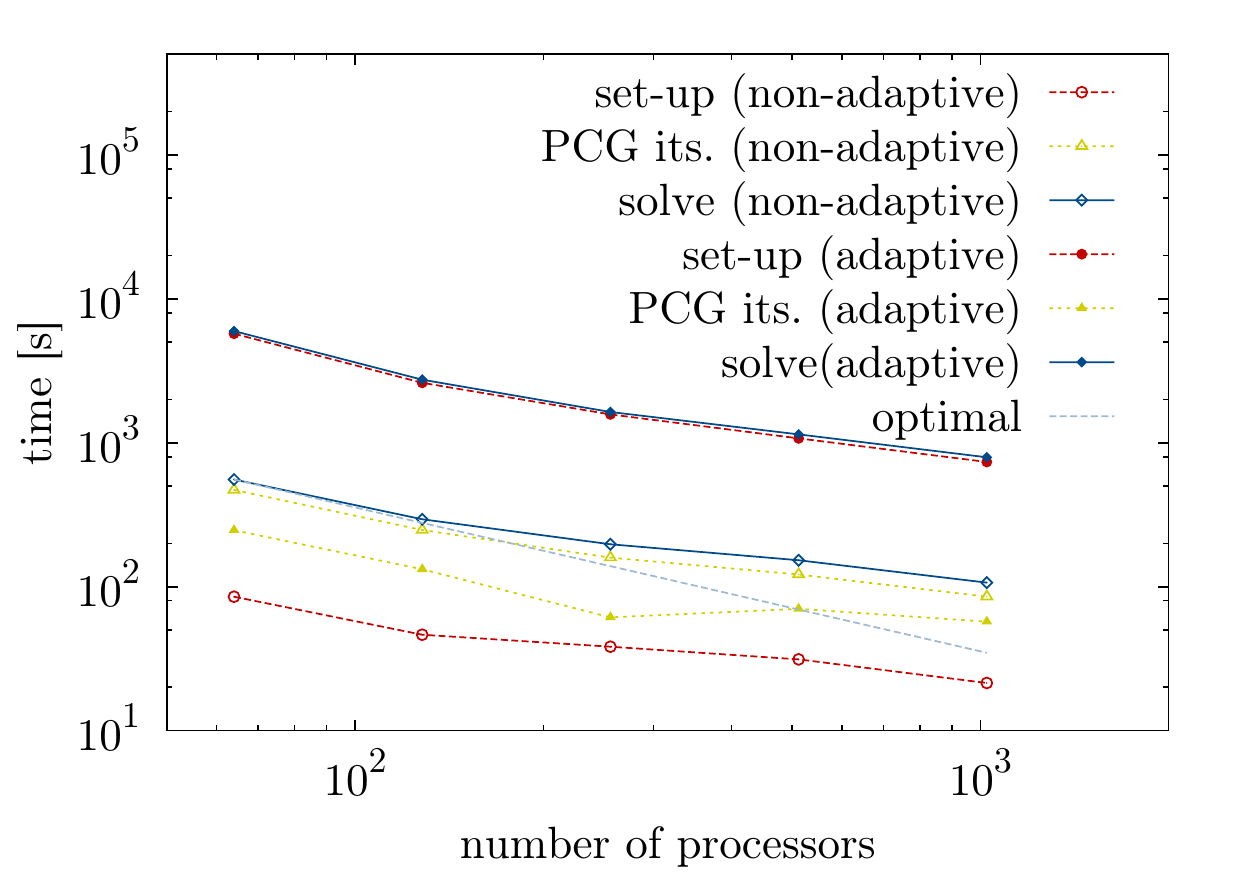}
\includegraphics[width=0.495\textwidth]{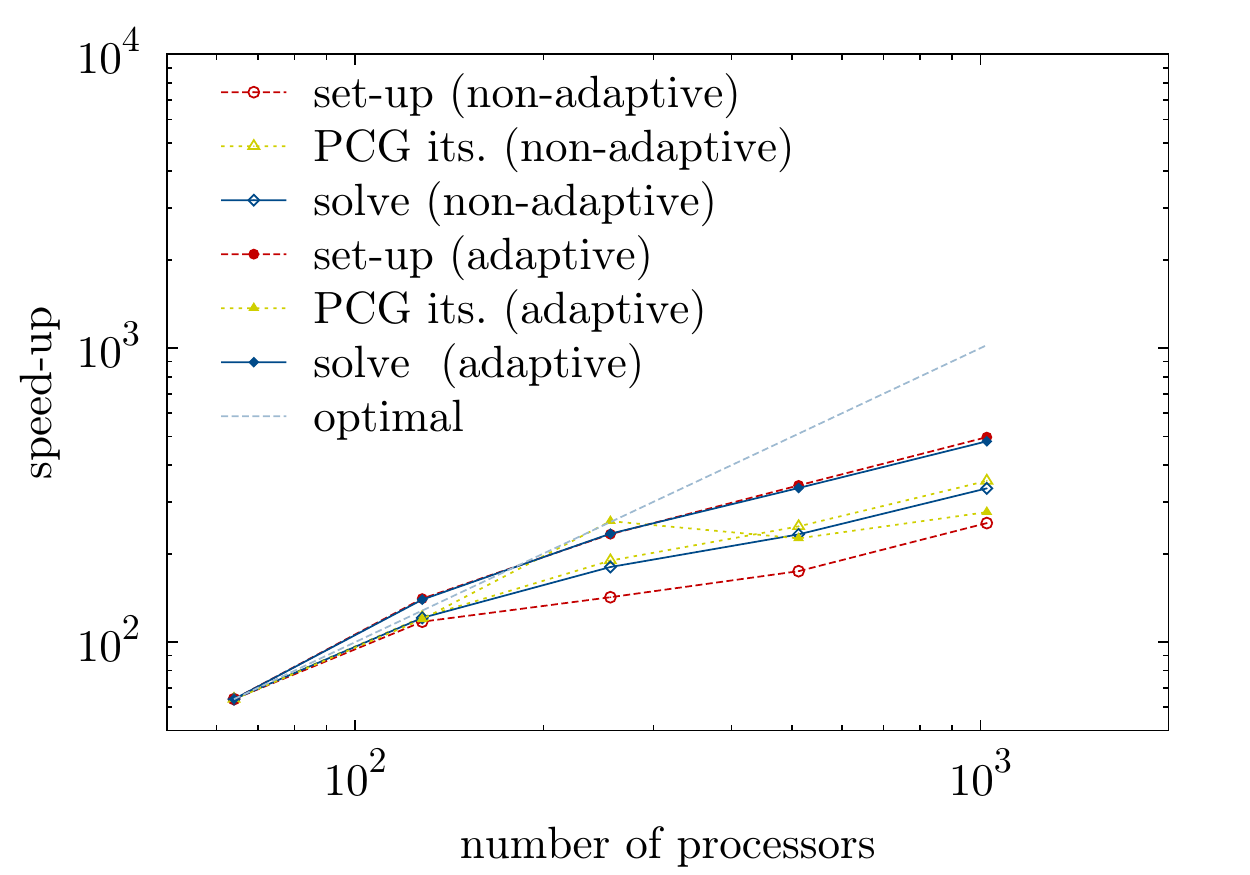}
\caption{\label{fig:timing_geoc1024_3level}
Strong scaling test for the geocomposite problem (Fig.~\ref{fig:geoc}) containing approx. 6 million degrees of freedom,
on the first level divided into 1024 subdomains with 5635 faces with arithmetic/adaptive constraints, 
and on the second level into 32 subdomains with 100 faces.
Computational time (left) and speed-up (right) separately for set-up and PCG phases, and their sum (`solve'),
comparison of \emph{non-adaptive} (214 its.) and \emph{adaptive} (113 its.) \emph{parallel 3-level BDDC}. 
}
\end{center}
\end{figure}

In Figs.~\ref{fig:timing_geoc1024_2level} and \ref{fig:timing_geoc1024_3level}, 
we report strong scaling of the parallel implementation on this problem for \emph{2-level} and \emph{3-level BDDC}, 
respectively.
The strong scaling is again nearly optimal. 

We can again see that while for the non-adaptive approach, most time is spent in the PCG iterations,
for the adaptive approach, the curve for the set-up phase is almost indistinguishable from the one for total solution time, 
and the set-up clearly dominates the solution.
For the two-level method, MUMPS was not able to solve the coarse problem on 1024 cores and so this value is omitted in 
Fig.~\ref{fig:timing_geoc1024_2level}.

\section{Conclusion}

\label{sec:conclusion}

We have presented the algorithm of \emph{Adaptive-Multilevel BDDC} method for three-dimensional problems, 
and its parallel implementation.
The algorithm represents the concluding step of combining ideas developed separately for adaptive selection of constraints in BDDC~
\cite{Mandel-2007-ASF,Mandel-2012-ABT,Sistek-2012-SPA,Sousedik-2008-CDD},
and for the multilevel extension of the BDDC method~\cite{Mandel-2007-OMB,Mandel-2008-MMB,Sistek-2012-PIM}.
The algorithm and its serial implementation was studied in the thesis \cite{Sousedik-2010-AMB-thesis}, 
and the serial algorithm for two dimensional problems also in \cite{Sousedik-2011-AMB}.

The \emph{Adaptive BDDC} method aims at numerically difficult problems, like those containing severe jumps in material coefficients
within the computational domain.
It recognizes troublesome parts of the interface 
by solving a generalized eigenvalue problem for each pair of adjacent subdomains which share a face.
By dominant eigenvalues the method detects where constraints need to be concentrated in order to improve the coarse space, thus reducing number of iterations. 
On the other hand, the \emph{Multilevel BDDC} aims at improving scalability of the BDDC method for very large numbers of subdomains, 
for which the coarse problem gets too large and/or fragmented to be solved by a parallel direct solver.
However, as theory suggests and experiments confirm, 
\emph{Multilevel BDDC} leads to an exponential growth of the condition number and the number of iterations.

The \emph{Adaptive-Multilevel BDDC} method
provides a kind of synergy of the adaptive and the multilevel approaches.
Our results confirm, 
that adaptively generated constraints are capable 
of reducing the rate of growth of the condition number with levels.
At the same time, the extension to three levels improved scalability of the adaptive 2-level approach, 
and for large problems and core counts, we have been able to obtain results we could not get by the 2-level method.

A convenient way of preconditioning the LOBPCG method based on 
components of BDDC was presented,
effectively converting the generalized eigenvalue problems to ordinary eigenproblems.
In our computations, this preconditioning led to large savings of number of LOBPCG iterations and corresponding computing time. 
However, to reduce the time necessary for solving the eigenproblems further, 
we have restricted the maximal number of adaptively generated constraints per face to ten in our parallel computations.
For this reason, 
the resulting performance of the adaptive method was not as optimal as the serial tests (with arbitrary number 
of constraints per face) suggested.

We have described a~parallel implementation of the algorithm available in our open-source library BDDCML.
The solver has been successfully applied to systems of equations with over 400 million unknowns solved on 32 thousand cores.
Presented results confirm, that both adaptive and non-adaptive implementations are reasonably scalable. 
However, presented computations have also revealed sub-optimal scaling especially on large numbers of cores,
and these results will provide a basis for further optimization of the solver.

We have presented results of two benchmark and two engineering problems of structural analysis.
On all problems, the adaptive selection of constraints led to reduced number of PCG iterations.
However, for most problems,
this fact did not lead to savings in computational time, and the cost of generating adaptive constraints was not compensated
by saved iterations.

It can be concluded that the adaptive method is not suitable for simple problems, 
where also non-adaptive (even multilevel) method would converge reasonably fast.
For problems with difficulties, the \emph{non-adaptive BDDC} method leads to large cost of PCG iterations, 
especially for using several levels.
On the contrary, the set-up phase with solution of local eigenproblems
mostly dominates the overall solution time for the \emph{Adaptive-Multilevel BDDC} method.
Which approach is finally advantageous depends on a particular problem.
Apart of the aspect of computational time, 
we have already encountered several problems for which the \emph{non-adaptive BDDC} method failed,
and which have been successfully solved by \emph{adaptive BDDC}.

\section*{Acknowledgement}
We would like to thank to Jaroslav Novotn\'{y}, Jan Le\v{s}tina, Radim Blaheta, and Ji\v{r}\'{\i} Star\'{y} 
for providing data of real engineering problems. 

This work was supported in part by National Science Foundation under grant 
\mbox{DMS-1216481}, by Czech Science Foundation under grant \mbox{GA \v{C}R 106/08/0403}, 
and by the Academy of Sciences of the Czech Republic through \mbox{RVO:67985840}.
B. Soused\'{\i}k acknowledges support from the DOE/ASCR and the NSF PetaApps award number 0904754.
J. {\v S}\'{\i}stek acknowledges the computing time on \emph{Hector} supercomputer provided by the PRACE-DECI initiative.
A part of the work was done at the University of Colorado Denver when B. Soused\'{\i}k was a graduate student 
and during visits of J. {\v S}\'{\i}stek, partly supported by the Czech-American Cooperation program of
the Ministry of Education, Youth and Sports of the Czech Republic under research project LH11004.


\end{document}